\documentclass[final,leqno,onefignum,onetabnum]{siamltex1213}

\usepackage{cite,extarrows,bm,amssymb,amsmath,enumerate,color,txfonts,mathrsfs,tipa,multirow,epsfig,indentfirst}%
\usepackage{float}

\def\be{\begin{equation}}
\def\en{\end{equation}}
\def\beq{\begin{eqnarray}}
\def\eq{\end{eqnarray}}
\def\beqx{\begin{eqnarray*}}
\def\eqx{\end{eqnarray*}}

\def\vv{\mbox{\sc v}}

\makeatletter
  \newcommand\figcaption{\def\@captype{figure}\caption}
  \newcommand\tabcaption{\def\@captype{table}\caption}
\makeatletter
\newtheorem{remark}{Remark}[section]

\title{Plane Wave Discontinuous Galerkin methods for the Helmholtz equation and Maxwell equations in Anisotropic Media
\thanks{The first author was supported by China NSF under the grant 11501529 and Qinddao
applied basic research project under grant 17-1-1-9-jch.
The second author was supported by the Natural Science Foundation of China G11571352.}}
\author{LONG YUAN
\thanks{ College of Mathematics and Systems Science, Shandong
University of Science and Technology,
Qingdao 266590, China (yuanlong@lsec.cc.ac.cn).}
 \and Qiya Hu
\thanks{ 1. LSEC, Institute of Computational Mathematics and
Scientic/Engineering Computing, Academy of Mathematics and Systems
Science, Chinese Academy of Sciences, Beijing 100190, China; 2. School of Mathematical Sciences, University of Chinese Academy
of Sciences, Beijing 100049, China
(hqy@lsec.cc.ac.cn).}
}

\begin{document}

\maketitle

\begin{abstract}
In this paper we are concerned with plane wave discontinuous Galerkin (PWDG) methods for Helmholtz equation and time-harmonic Maxwell equations
in three-dimensional anisotropic media, for which the coefficients of the equations are matrices instead of numbers.
We first define novel plane wave basis functions based on rigorous choices of scaling transformations and coordinate transformations. Then we
derive the error estimates of the resulting approximate solutions
with respect to the condition number of the coefficient matrices, under a new assumption on the shape regularity of polyhedral meshes.
Numerical results verify the validity of the theoretical results,
and indicate that the approximate solutions generated by the proposed PWDG method
\textcolor{red}{possess high accuracy.}
\end{abstract}



\pagestyle{myheadings} \thispagestyle{plain} \markboth{ LONG YUAN and QIYA HU}{PWDG Methods in Anisotropic Media}

\section{Introduction}

The plane wave method, which is based on the Trefftz approximation space made of plane wave basis functions,
 was first introduced to discretization of homogeneous Helmholtz equation and was then extended to discretization of
homogeneous time-harmonic Maxwell's equations and elastic wave equations. Various examples of the plane wave methods has been systematically surveyed in \cite{hmp},
 for example, the ultra weak variational formulation (UWVF) \cite{refOC,ref11,hmm,hm,Gabard}, the plane wave discontinuous Galerkin (PWDG) method \cite{ref20,ref21,HMP_3M,HMP_apnum,pwdg,hmp,HMP_FOUND},
 the plane wave least-squares (PWLS) method \cite{ref12,hy,hy2,yh2,peng}. The plane wave least-squares
 combined with local spectral finite element (PWLS-LSFE) method was recently proposed in \cite{hy3} for nonhomogeneous Helmholtz equation
and time-harmonic Maxwell equations. The plane wave method has
 an important advantage over Lagrange finite elements for discretization of the Helmholtz equation
 and time-harmonic Maxwell equations \cite{ref21,hmp,hmm}: to achieve the same accuracy, \textcolor{red}{relatively fewer degrees of freedom} are enough in the plane wave-type methods owing to
 the particular choice of the basis functions that \textcolor{red}{(possibly approximately)} satisfy the considered PDE without boundary conditions.

In \cite{hm} the UWVF method was first extended to discretization of homogeneous
Maxwell's equations in anisotropic media. The studies in \cite{hm} were
devoted to approximating the Robin-type trace of the electric and magnetic fields in an anisotropic medium, and focus on the numerical tests and convergence analysis in TM mode scattering, which can result in a Helmholtz equation in two dimensions with an anisotropic coefficient. \textcolor{red}{ In addition, the work \cite{hm,yuan,yuan2} proposed the anisotropic plane wave basis functions, based on
a coordinate transformation, for two-dimensional anisotropic Helmholtz equations.
}

It was pointed out in \cite[p.351]{hm} that, for three-dimensional anisotropic time-harmonic Maxwell's system, {\it almost all theoretical questions
related to the 3D UWVF approach to anisotropic media are still open:
in particular, the relevant approximation properties of sums of
anisotropic plane waves are not known.} \textcolor{red}{Meanwhile, the three-dimensional Maxwell equations in anisotropic media with positive definite matrices play an important role in practical physical applications (see \cite[p. 5-6]{monk}).}
 Recently, the PWDG method was applied to
discretization of three-dimensional anisotropic time-harmonic Maxwell's equations with diagonal matrix coefficients \cite{yuan_aniso},
in which rough error estimates of the resulting approximate solutions were derived.


In this paper, we \textcolor{red}{first} study the PWDG method for three-dimensional Helmholtz equation and time-harmonic Maxwell equations in more general anisotropic media,
where the coefficient of the equations is a positive definite matrix instead of diagonal matrix. In order to deal with such complicated models and build better
convergence results, we have to carefully define plane wave basis functions by rigorous choices of the scaling transformations and the coordinate transformations. Besides, we propose a new assumption on the triangulation: the transformed triangulation
$\hat{{\cal T}}_{\hat{h}}$ rather than the physical triangulation ${\cal T}_h$ is shape regular. Own to such an assumption, we can
verify that the defined transformations have the desired stability estimates on the condition number of the anisotropic coefficient matrix, and further
prove that the resulting approximate solutions have the desired convergence order with respect to the condition number of the coefficient matrix.
Numerical experiments indicate that the approximate solutions generated
by the proposed PWDG method \textcolor{red}{possess high accuracy}, and verify the validity of the theoretical results.

\textcolor{red}{In order to clarify the novelty of the present paper, we would like to emphasize the differences with the existing works:}

\textcolor{red}{1) the models are different: in the present paper, we consider three-dimensional time-harmonic Maxwell's equations with general coefficients that are (constant) positive definite matrices,
but in [26] three-dimensional time-harmonic Maxwell's equations with diagonal matrix coefficients were considered;
in [16], three-dimensional time-harmonic Maxwell's equations that can be transformed into two-dimensional Helmholtz equations were studied;
in [25, 24] two-dimensional Helmholtz equations with coefficients
being positive definite matrices were studied.
}

\textcolor{red}{2) the methods are different: in the present paper, we employ a new transformation $\hat{\bf x} = \Lambda^{-\frac{1}{2}}P {\bf x}$ that is different from the old transformation
$\hat{\bf x} = A^{-\frac{1}{2}} {\bf x}$ that used in [16, 26, 25, 24].
Moreover, we propose a new assumption on the triangulation such that the transformed triangulation $\hat{{\cal T}}_{\hat{h}}$ rather than the physical triangulation ${\cal T}_h$ is shape regular.}

\textcolor{red}{3) the results are different: in the present paper, we obtained obviously better error estimates than the existing works, thanks to the proposed transformation.}

The paper is organized as follows. In Section 2, we present PWDG method for the three-dimensional anisotropic Helmholtz equation and analyze
the convergence of the proposed PWDG method.
In Section 3, we extend the results obtained in Section 2 to the three-dimensional anisotropic time-harmonic Maxwell equations.
Finally, we report some numerical results to confirm the effectiveness of the proposed method.

\section{PWDG method for three-dimensional anisotropic Helmholtz equation}
To our knowledge, there seems no work on plane wave method for three-dimensional anisotropic Helmholtz equation in literature.
In this section we first introduce PWDG method for three-dimensional anisotropic Helmholtz equation.


\subsection{The model and its variational formula}
Consider the three-dimensional (3D) anisotropic Helmholtz equation of acoustic wave field $u$ (refer to \cite{Imbert}):
\be \label{3danisohelm}
 -\nabla\cdot A\nabla u - \omega^2 u = 0  \quad \text{in} \quad \Omega
\en
with the Robin-type boundary condition
\be \label{robin}
{\bf n}\cdot A\nabla u + i\omega u = g \quad \text{on} \quad\gamma=\partial\Omega.
\en
Here $\Omega$ is a bounded domain in three dimensions, ${\bf n}$ denote the unit outer normal vector to the boundary $\partial \Omega$; $A$ is a real-valued $3 \times 3$ matrix; \(\omega>0\)
is the temporal frequency of the field and \({ g}\in L^2(\partial\Omega)\).
We assume that the matrix $A$ is symmetric and positive definite.

We will derive a variational formula of the considered model based on a partition of the solution domain $\Omega$ (refer to \cite{ref21}).

For convenience, assume that $\Omega$ is a polyhedron. Let \(\Omega\) be divided into a union of some
elements in the sense that
$$ \overline{\Omega}=\bigcup_{k=1}^N\overline{\Omega}_k,\quad
\Omega_l\bigcap\Omega_j=\emptyset
 \quad\text{ for }l\not=j,$$
where each $\Omega_k$ is a polyhedron. Let \( {\cal T}_h\) denote
the partition comprised of the elements \(\{\Omega_k\}\), where
\(h\) denotes the diameter of the maximal element in \(\{\Omega_k\}\). Define $$
\Gamma_{lj}=\partial\Omega_l\bigcap\partial\Omega_j, \quad\text{for
 }l\not=j$$
 and
$$ \gamma_k=\overline{\Omega}_k\bigcap\partial\Omega
 \quad(k=1,\ldots,N),~~
 \gamma=\bigcup^N_{k=1}\gamma_k.$$
We denote by $\mathcal{F}_h = \bigcup\limits_k\partial\Omega_k$ the
skeleton of the mesh, and set $\mathcal{F}_h^{\text{B}}=\mathcal{F}_h\bigcap\partial\Omega$ and
$\mathcal{F}_h^{\text{I}}= \mathcal{F}_h \backslash
\mathcal{F}_h^{\text{B}}$.
Let $u$ and ${\bm \sigma}$ be a piecewise smooth function and vector field on ${\cal T}_h$, respectively.
On $\partial\Omega_l\bigcap\partial\Omega_j$, we define
\begin{eqnarray}
& \text{the averages:} ~~\{\{ u \}\}:= \frac{u_l+u_j}{2}, ~~\{\{
{\bm \sigma}\}\} := \frac{{\bm \sigma}_l + {\bm \sigma}_j}{2},
\cr &
\text{the jumps:} ~\llbracket {u} \rrbracket_N := u_l{\bf n}_l + u_j{\bf n}_j,
~\llbracket {\bm \sigma} \rrbracket_T:= {\bf n}_l\times {\bm \sigma}_l + {\bf n}_j\times {\bm \sigma}_j,
\cr &
~\llbracket {\bm \sigma} \rrbracket_N = {\bf n}_l\cdot {\bm \sigma}_l + {\bf n}_j\cdot {\bm \sigma}_j,
\end{eqnarray}
where ${\bf n}$ denotes the unit outer normal vector on the boundary of each element \(\Omega_k\).

Define the broken Sobolev space
\be
H^s({\cal T}_h)=\{ v\in L^2(\Omega):~ v|_{\Omega_k} \in H^s(\Omega_k)~~ \text{for} ~~\forall \Omega_k\in {\cal T}_h \}.
\en
As usual, we assume that each entry in $A$ is a constant. Let ${V}({\cal T}_h)$ be the piecewise Trefftz space defined on ${\cal T}_h$ by
\be
{V}({\cal T}_h) = \{v\in H^2({\cal T}_h);~-\nabla\cdot A\nabla v - \omega^2 v =0 ~~\text{in ~~each}~~\Omega_k \in {\cal T}_h \}.
\en

Let $\alpha$ and $\beta$ be two positive numbers, and let $\delta\in (0,\frac{1}{2}]$.
Define the sesquilinear form $A_h(\cdot, \cdot)$ by
\begin{eqnarray}\label{helmpwdg9}
& \mathcal{A}_h(u, v) = \int_{\mathcal{F}_h^{\text{I}}} \{\{ u\}\}~
\overline{\llbracket A \nabla_h v \rrbracket_N}~dS
+i\omega^{-1}\int_{\mathcal{F}_h^{\text{I}}}
\beta \llbracket A \nabla_h u \rrbracket_N \cdot
\overline{\llbracket A \nabla_h v \rrbracket_N}~dS
  \cr &
- \int_{\mathcal{F}_h^{\text{I}}} \{\{ A \nabla_h u \}\}\cdot
\overline{\llbracket v \rrbracket_N}~dS
+ i\omega\int_{\mathcal{F}_h^{\text{I}}}\alpha
\llbracket  u \rrbracket_N \cdot \overline{\llbracket v \rrbracket_N}~dS
\cr & +
\int_{\mathcal{F}_h^{\text{B}}}(1-\delta)u ~\overline{( A \nabla_h v \cdot{\bf n} )}~dS
+i\omega^{-1}\int_{\mathcal{F}_h^{\text{B}}} \delta
( A \nabla_h u \cdot {\bf n}) ~ \overline{( A \nabla_h v \cdot{\bf n})}~dS
\cr &
 -\int_{\mathcal{F}_h^{\text{B}}} \delta ( A \nabla_h u \cdot{\bf n}) ~\overline{v}~dS
 + i\omega\int_{\mathcal{F}_h^{\text{B}}} (1-\delta)u~ \overline{v}~dS,
 \quad \forall u,v \in V({\cal T}_h).
\end{eqnarray}
and \begin{eqnarray} \label{helmpwdg10} &
\ell_h(v)=
i\omega^{-1}\int_{\mathcal{F}_h^{\text{B}}}
\delta g ~ \overline{ ( A \nabla_h v \cdot{\bf n}) }~dS
+ \int_{\mathcal{F}_h^{\text{B}}} (1-\delta)g~\overline{v}~dS, \quad \forall v\in V({\cal T}_h).
 \end{eqnarray}

\textcolor{red}{Then, for a given $g$, the variational problem associated with (\ref{3danisohelm})-(\ref{robin}) can be expressed as follows (see \cite[Section 2]{ref21}): Find $u \in V({\cal T}_h)$ such that}
\be \label{helmanisohelmvaria}
 \mathcal{A}_h(u, v) = \ell_h(v), \quad \forall v\in V({\cal T}_h).
\en

We endow the space ${V}({\cal T}_h)$ with the norm
\begin{eqnarray} \label{helmfnorm}
&|||v|||_{\mathcal{F}_h}^2~:=
\omega^{-1}||\beta^{1/2}\llbracket A \nabla_h v\rrbracket_N||_{0,\mathcal{F}_h^I}^2
+ \omega||\alpha^{1/2}\llbracket v\rrbracket_N||_{0,\mathcal{F}_h^I}^2
\cr & + \omega^{-1}||\delta^{1/2} A \nabla_h v\cdot {\bf
n}||_{0,\mathcal{F}_h^B}^2 +
\omega||(1-\delta)^{1/2}v||_{0,\mathcal{F}_h^B}^2
\end{eqnarray}
and the augmented norm
\begin{eqnarray} \label{helmfplusnorm}
& |||v|||_{\mathcal{F}_h^+}^2~:= |||v|||_{\mathcal{F}_h}^2 +
\omega||\beta^{-1/2}\{\{ v\}\}||_{0,\mathcal{F}_h^I}^2 \cr & +
\omega^{-1}||\alpha^{-1/2}\{\{ A \nabla_h v
\}\}||_{0,\mathcal{F}_h^I}^2 +
\omega||\delta^{-1/2}v||_{0,\mathcal{F}_h^B}^2.
\end{eqnarray}

As in \cite{ref21}, we can show the following existence, uniqueness and continuity
of solution of the above variational problem.
\begin{lemma} \label{helmexistetc}
There exists a unique $u$ solution to (\ref{helmanisohelmvaria}); moreover, we have
 \begin{eqnarray} \label{helmuniquebound}
& -\text{Im}[\mathcal{A}_h({ w},{ w})] =
\big|\big|\big|{ w}\big|\big|\big|^2_{\mathcal{F}_h}  \cr &
\text{and} ~~\big|\mathcal{A}_h({ w}, \xi)\big| \leq C
\big|\big|\big|{ w}\big|\big|\big|_{\mathcal{F}_h^+}~
\big|\big|\big|{ \xi}\big|\big|\big|_{\mathcal{F}_h}, ~~\forall ~{
w},{ \xi}\in V({\cal T}_h).
\end{eqnarray}
\end{lemma}

\subsection{Discretization for the 3D anisotropic Helmholtz equation}
\textcolor{red}{
Since $A$ is positive definite matrix, there exists an orthogonal matrix $P$ and a diagonal positive definite matrix $\Lambda = \text{diag} ( \lambda_{\text{min}},\lambda_{\text{mid}},\lambda_{\text{max}} )$ such that
$A= P^T \Lambda P$, where $\lambda_{\text{min}}\leq \lambda_{\text{mid}} \leq \lambda_{\text{max}}$ are constant and the superscript $T$ denotes matrix transposition.
}
Of course, we can assume that $\text{det}(P)=1$. It is clear that $A^{\frac{1}{2}}= P^T \Lambda^{\frac{1}{2}}P$.
Define a coordinate transformation:
\be\label{tran1}
 \hat{\bf x} = \Lambda^{-\frac{1}{2}}P{\bf x} \xlongequal{\Delta} S{\bf x}, \quad \hat{\bf x} =(\hat{x} ~ \hat{y} ~ \hat{z})^T, \quad S=\Lambda^{-\frac{1}{2}}P,\quad {\bf x}=(x~y~z)^T.
\en
For convenience, we use $p_1,p_2$ and $p_3$ to denote three column vectors of $P$, and use $q_1,q_2,q_3$ to denote the three row vectors of $P$. Then each of these vectors is a three dimensional
unit vector. Moreover, $p_1,p_2$ and $p_3$ (and $q_1,q_2,q_3$) are orthogonal each other.
Set $u({\bf x}) = u( S^{-1}\hat {\bf x}) \xlongequal{\Delta} \hat u(\hat {\bf x})$.

Let $\hat{\Omega}$ and $\hat{\Omega}_k$ denote the images of $\Omega$ and $\Omega_k$ under the coordinate transformation (\ref{tran1}), respectively.
Since the map $S$ is linear, the transformed domain $\hat\Omega$ and elements $\hat{\Omega}_k$ are also
polyhedrons. We use $\hat{\bf n}$ to denote the unit outer normal vector on the boundary of each element \(\hat\Omega_k\), and $\hat{{\cal T}}_{\hat{h}}$ to denote the
partition comprised by the elements $\{\hat{\Omega}_k\}$, where $\hat{h}$ is the maximal diameter of the
transformed elements $\{\hat{\Omega}_k\}$. Set
$\hat{\mathcal{F}}_{\hat h} = \bigcup\limits_k\partial\hat\Omega_k$,
$\hat{\mathcal{F}}_{\hat h}^{\text{B}}=\partial\hat{\Omega}$ and
$\hat{\mathcal{F}}_{\hat h}^{\text{I}}= \hat{\mathcal{F}}_{\hat h}
\backslash \hat{\mathcal{F}}_{\hat h}^{\text{B}}$.

We denote by $\nabla_h$ and $\hat\nabla_{\hat h}$ the element application of the gradient operator $\nabla=(\frac{\partial}{\partial x} ~\frac{\partial}{\partial y} ~\frac{\partial}{\partial z} )^T$ and $\hat\nabla=(\frac{\partial}{\partial \hat x} ~\frac{\partial}{\partial \hat y} ~\frac{\partial}{\partial \hat z} )^T$, respectively.
Define the Laplace operator $\hat\triangle$ on $\hat{\Omega}$ by $\hat\triangle=\frac{\partial^2}{\partial\hat x^2}+\frac{\partial^2}{\partial\hat y^2}+\frac{\partial^2}{\partial\hat z^2}$.

\textcolor{red}{By direct calculation,} we have
$$A\nabla u = P^T \Lambda^{\frac{1}{2}} \nabla \hat u = S^{-1} \hat\nabla \hat u, \quad\text{and} \quad \nabla\cdot A\nabla u = \hat\nabla\cdot (\hat\nabla \hat u ) \xlongequal{\Delta} \hat\triangle \hat u.$$
Thus, the anisotropic Helmholtz equation (\ref{3danisohelm}) is transformed into the isotropic Helmholtz equation
\be \label{helm8}
 \hat\triangle \hat u + \omega^2 \hat u = 0 \quad \text{in} \quad \hat\Omega.
 \en

 Conversely, if $\hat u$ satisfies the isotropic Helmholtz equation (\ref{helm8}),
 the function $u({\bf x}) = \hat u( S {\bf x})$ satisfies the original anisotropic Helmholtz equation (\ref{3danisohelm}).

In order to define suitable anisotropic plane wave basis functions, we first define
plane wave basis functions $\{y_{kl}\}$ satisfying the isotropic Helmholtz equation (\ref{helm8}) on $\hat\Omega_k$ as follows.
\begin{eqnarray}
\left\{\begin{array}{ll}  y_{kl}(\hat{\bf x}) =
e^{i\omega(\hat{\bf x}\cdot \boldsymbol{d_l})},~~\hat{\bf x}\in\hat\Omega_k,\\
\boldsymbol{d_l}\cdot \boldsymbol{d_l}=1,\\
l\neq s \rightarrow
\boldsymbol{d_{l}}\neq\boldsymbol{d_{s}},
\end{array}\right.
\label{helm10}
\end{eqnarray}
where $\boldsymbol{d_l}~(l=1,\cdots,p)$ are unit wave
propagation directions, and can be determined by the codes in \cite{refsite}.
Choose the number $p$ of plane wave propagation
directions as $p=(m+1)^2$, where $m$ is a positive integer.

Then the anisotropic plane wave basis functions of $V_p({\mathcal T}_h)$ can be defined as
\begin{eqnarray}
\quad\quad u_{kl}({\bf x})=\left\{\begin{array}{ll} y_{kl}(S{\bf x}),
~~{\bf x}\in \Omega_k,\\
 0,~~{\bf x}\in \Omega_j~~\mbox{satisfying}~~j\neq k
\end{array}\right.~~(k,j=1,\cdots,N;~l=1,\cdots,p).
\label{helm11}
\end{eqnarray}
Thus the space $V({\mathcal T}_h)$ is discretized by the subspace
\begin{equation}
V_p({\mathcal T}_h)=\text{span}\bigg\{u_{kl}:~k=1,\cdots,N;~l=1,\cdots,p \bigg\}. \label{helm12}
\end{equation}

Furthermore, we obtain the discretized version of the continuous variational problem (\ref{helmanisohelmvaria}):
Find $u_h \in V_p({\cal T}_h)$ such that
\be \label{helmdisvaria}
 \mathcal{A}_h(u_h, v_h) = \ell_h(v_h), \quad \forall v_h\in V_p({\cal T}_h).
\en

\subsection{Error estimates of the approximate solutions}
In this subsection, we are devoted to the analysis of convergence of the plane wave approximation $u_h$.

\subsubsection{The required partition}
In order to derive the desired error estimates of the approximate solutions, we require that the partition must satisfy some assumptions. In this part we introduce a kind of particular triangulation
such that these assumptions can be met.

\begin{figure}[H]
\begin{center}
\begin{tabular}{c}
 \epsfxsize=0.7\textwidth\epsffile{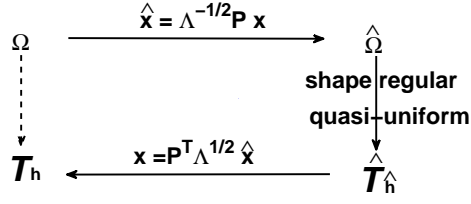}\\
\end{tabular}
\end{center}
 \caption{Mesh generation. } \label{mesh}
\end{figure}

\textcolor{red}{We adopt a non-regular triangulation} \(\mathcal {T}_h\) for the three-dimensional \textcolor{red}{convex} domain \(\Omega\) as follows (see Figure \ref{mesh}).

\noindent {\bf Mesh Generation Algorithm:}

{\it Step 1.}  Determine the image domain $\hat\Omega$ of $\Omega$ under the coordinate transformation (\ref{tran1}).

{\it Step 2.}  Decompose $\hat\Omega$ into polyhedron elements $\{\hat{\Omega}_k\}$ such that $\hat{{\cal T}}_{\hat{h}}$ is shape regular and quasi-uniform
in the usual manner. \textcolor{red}{Besides, we assume that each element $\hat{\Omega}_k$ of $\hat{{\cal T}}_{\hat{h}}$ is a convex Lipschitz domain.
}

{\it Step 3.}  Determine the triangulation \( {\cal T}_h\) of $\Omega$ by using the inverse transformation of (\ref{tran1}) acting on the elements of $\hat{{\cal T}}_{\hat{h}}$.

\begin{remark} Notice that the coordinate transformation (\ref{tran1}) consists of the orthogonal transformation $P$ and the scaling transformation $\Lambda^{-\frac{1}{2}}$.
It is known that an orthogonal matrix preserves the length of any vector and the angle between two vectors, then the orthogonal transformation preserves the shape of any polyhedron.
Thus, it is cheap to determine the image domain $\hat\Omega$ and the triangulation \( {\cal T}_h\) by the coordinate transformation (\ref{tran1}) and its inverse transformation, respectively, since we need only to compute the coordinates of the vertices of the polyhedron $\hat\Omega$ and the elements $\Omega_k$.
\end{remark}

For a $3 \times 3$ matrix $B$, we define
\be
||B||=||B||_2=\mathop{\text{max}~~}\limits_{{\bf 0} \neq {\bf x}\in R^3}\frac{||B{\bf x}||_2}{||{\bf x}||_2},
\en
where ${\bf x}=(x_1, x_2, x_3)^T$ and $||{\bf x}||_2=(x_1^2+x_2^2+x_3^2)^{\frac{1}{2}}$.
Then it is well known that $||A||=||\Lambda||=||\Lambda^{\frac{1}{2}}||^2$.
For the simplicity of notation, let $\rho$ denote the condition number $\text{cond}(A)$ of the anisotropic matrix $A$.
Then $\rho = \text{cond}(\Lambda) = \text{cond}^2(\Lambda^{\frac{1}{2}})= \text{cond}^2(S)$.

The following Lemma gives an important geometric property of the triangulations.

\begin{lemma} For the proposed triangulation, we have
\be \label{geometric100}
c_0 ||\Lambda^{\frac{1}{2}}||^{-1}h \leq \hat{h} \leq C_0||\Lambda^{\frac{1}{2}}||^{-1}h,
\en
where $c_0$ and $C_0$ denote two constants independent of $\omega, \rho, h, p$.
\end{lemma}

{\it Proof}. Let us first recall the definition of the shape regularity and quasi-uniformity assumptions.
There exists a constant $C$ independent of $\hat\Omega_k$ and $\hat{{\cal T}}_{\hat{h}}$ such that for all $\hat\Omega_k \in \hat{{\cal T}}_{\hat{h}}$ and all $\hat{{\cal T}}_{\hat{h}}$,
\be  \label{hatassumpuni}
\frac{\hat h_k}{\hat \rho_k} \leq C \quad\mbox{and}\quad
\frac{\hat h}{\hat h_k} \leq C,
 \en
where $\hat h_k$ denotes the diameter of the minimum sphere containing the polyhedron $\hat\Omega_k$,
$\hat \rho_k$ denotes the diameter of the maximum sphere contained by the polyhedron $\hat\Omega_k$.

For convenience, we define two sub-transformation $\tilde{\bf x}=P{\bf x}$ and $\hat{\bf x}= \Lambda^{-\frac{1}{2}}\tilde{\bf x}$ to achieve the coordinate transformation (\ref{tran1}).
At first we consider the transformation
\be  \label{firtran}
\hat{\bf x} = (\hat x ~ \hat y ~ \hat z)^T= \Lambda^{-\frac{1}{2}} \tilde{\bf x} = ( \tilde x~ \tilde y ~ \tilde z )^T.
\en

For each element $\hat{\Omega}_k \in \hat{{\cal T}}_{\hat{h}}$, we use $\hat{\vv}_k^{(i)}$ $(i=1,2,\cdots)$ to denote all the vertices of the polyhedron $\hat{\Omega}_k$,
and assume that $\hat{\Omega}_k$ and $\hat{\vv}_k^{(i)}$ are transformed to $\tilde{\Omega}_k$ and $\tilde{\vv}_k^{(i)}$
under the inverse transformation of (\ref{firtran}), respectively.
We denote by $\tilde h_k$ the diameter of the minimum sphere containing the polyhedron $\tilde{\Omega}_k$,
and by $\tilde \rho_k$ the diameter of the maximum sphere contained by the polyhedron $\tilde{\Omega}_k$,
and set $\tilde h = \max\limits_k \tilde h_k$.

Define ${\Delta \hat x}_k,{\Delta \hat y}_k$ and ${\Delta \hat z}_k$ as
\be\label{axixlent}
{\Delta \hat x}_k = \max\limits_{i,j} \bigg(~\overrightarrow{\hat{\vv}_k^{(i)}\hat{\vv}_k^{(j)}}~\bigg)_{\hat x}, \quad
{\Delta \hat y}_k = \max\limits_{i,j} \bigg(~\overrightarrow{\hat{\vv}_k^{(i)}\hat{\vv}_k^{(j)}}~\bigg)_{\hat y}, \quad
{\Delta \hat z}_k = \max\limits_{i,j} \bigg(~\overrightarrow{\hat{\vv}_k^{(i)}\hat{\vv}_k^{(j)}}~\bigg)_{\hat z}.
\en
Here $\bigg(~\overrightarrow{\hat{\bf a}}~\bigg)_{\hat x} = a_1,~\bigg(~\overrightarrow{\hat{\bf a}}~\bigg)_{\hat y} = a_2$
 and $\bigg(~\overrightarrow{\hat{\bf a}}~\bigg)_{\hat z} = a_3$ for a vector $\overrightarrow{\hat{\bf a}} = (a_1 ~a_2~ a_3)^T$.
It follows by (\ref{hatassumpuni}) that
\be
\label{axixlent2}
 c \hat h_k \leq \hat \rho_k \leq {\Delta \hat x}_k \leq \hat h_k, \quad
 c \hat h_k \leq \hat \rho_k \leq {\Delta \hat y}_k \leq \hat h_k, \quad
 c \hat h_k \leq \hat \rho_k \leq {\Delta \hat z}_k \leq \hat h_k
\en
with $c = \frac{1}{C}$. Moreover, by the coordinate transformation (\ref{firtran}), we have
\be
\label{axixlent3}
{\Delta \hat x}_k= \frac{1}{\sqrt{\lambda_{\text{min}}}}{\Delta \tilde x}_k, \quad
{\Delta \hat y}_k= \frac{1}{\sqrt{\lambda_{\text{mid}}}}{\Delta \tilde y}_k, \quad
{\Delta \hat z}_k= \frac{1}{\sqrt{\lambda_{\text{max}}}}{\Delta \tilde z}_k,
\en
where
\be
\label{axixlent4}
{\Delta \tilde  x}_k= \max\limits_{i,j} \bigg(~\overrightarrow{\tilde{\vv}_k^{(i)} \tilde{\vv}_k^{(j)}}~\bigg)_{\tilde x}, \quad
{\Delta \tilde  y}_k= \max\limits_{i,j} \bigg(~\overrightarrow{\tilde{\vv}_k^{(i)} \tilde{\vv}_k^{(j)}}~\bigg)_{\tilde y}, \quad
{\Delta \tilde  z}_k= \max\limits_{i,j} \bigg(~\overrightarrow{\tilde{\vv}_k^{(i)} \tilde{\vv}_k^{(j)}}~\bigg)_{\tilde z}.
\en
\textcolor{red}{Combining} (\ref{axixlent2}) and (\ref{axixlent3}), yields
\be
c_1 ||\Lambda^{\frac{1}{2}}|| \hat h_k ~\leq ~\tilde h_k ~\leq ~C_1 ||\Lambda^{\frac{1}{2}}|| \hat h_k.
\en
Then, by (\ref{hatassumpuni}), we get
\be
c_2 ||\Lambda^{\frac{1}{2}}|| \hat h ~\leq ~\tilde h ~\leq ~C_2 ||\Lambda^{\frac{1}{2}}|| \hat h.
\en
Here $c_i$ and $C_i~(i=1,2)$ denote different constants independent of $\omega, \rho, h, p$ and the triangulation ${\mathcal T}_h$.

Finally, since the orthogonal transformation $\tilde{\bf x}=P{\bf x}$ does not change the length of vectors and the angle between vectors, we obtain the desired result.

$\Box$

\textcolor{red}{
The next Lemma gives a relation between the areas of two bounded planes based on the coordinate transformation  (\ref{tran1}).
}

\textcolor{red}{\begin{lemma} For the proposed triangulation, denote by $\Gamma$ a general bounded plane which belongs to  $\mathcal{F}_h$, and by $\hat\Gamma$ the corresponding plane belonging to $\hat{\mathcal{F}}_{\hat h}$. Then
we have
\be \label{arearela}
\frac{|\Gamma|}{|\hat{\Gamma}|}  \leq \lambda_{\text{mid}}^{\frac{1}{2}} \lambda_{\text{max}}^{\frac{1}{2}} ,
\en
where $|f|$ denotes the area of a bounded plane $f$ in the three-dimensional space.
\end{lemma}}

{\it \bf Proof}.  In order to achieve the coordinate transformation (\ref{tran1}), we also define two sub-transformation $\tilde{\bf x}=P{\bf x}$ and $\hat{\bf x}= \Lambda^{-\frac{1}{2}}\tilde{\bf x}$. Since the orthogonal transformation $\tilde{\bf x}=P{\bf x}$ does not change the length of vectors and the angle between vectors, we only need to prove the desired result under $ \hat{\bf x} = (\hat x ~ \hat y ~ \hat z)^T= \Lambda^{-\frac{1}{2}} \tilde{\bf x} = \Lambda^{-\frac{1}{2}} ( \tilde x~ \tilde y ~ \tilde z )^T$, namely,
\be \label{tr}
 \tilde{\bf x}= \Lambda^{\frac{1}{2}}\hat{\bf x}.
 \en

Without loss of generality, consider $\Gamma: \tilde z = a \tilde x + b \tilde y +c, ~ (\tilde x, \tilde y) \in D$, where $a,b,c$ are constants, $D\subset R^2$ is a bounded domain.  By (\ref{tr}), we get $\hat \Gamma: \hat z = a (\frac{\lambda_{\text{min}}}{\lambda_{\text{max}}})^{\frac{1}{2}} \hat x + b (\frac{\lambda_{\text{mid}}}{\lambda_{\text{max}}})^{\frac{1}{2}} \hat y + c (\frac{1}{\lambda_{\text{max}}})^{\frac{1}{2}} , ~(\hat x,\hat y) \in \hat D$. Obviously, $\frac{|D|}{|\hat D|} = det(\text{diag}(\lambda_{\text{min}}^{\frac{1}{2}},\lambda_{\text{mid}}^{\frac{1}{2}}))= \lambda_{\text{min}}^{\frac{1}{2}}\lambda_{\text{mid}}^{\frac{1}{2}}$.

Thus,
\be  \label{areasq1}
|\Gamma| = \iint_{D} \sqrt{1+\tilde z_{\tilde x}^2+ \tilde z_{\tilde y}^2 }~\mathrm{d}{\tilde x}\mathrm{d}{\tilde y} = \sqrt{1+a^2+ b^2 }~ |D|
\en
and
\be \label{areasq2}
|\hat\Gamma| = \iint_{\hat D} \sqrt{1+\hat  z_{\hat x}^2+ \hat z_{\hat y}^2 }~\mathrm{d}{\hat x}\mathrm{d}{\hat y} = \sqrt{1+a^2 \frac{\lambda_{\text{min}}}{\lambda_{\text{max}}}+ b^2\frac{\lambda_{\text{mid}}}{\lambda_{\text{max}}} }~ |\hat D|.
\en
 Combining the above equalities, we get
 \be
 \frac{|\Gamma|}{|\hat\Gamma|} = \lambda_{\text{min}}^{\frac{1}{2}}\lambda_{\text{mid}}^{\frac{1}{2}}\lambda_{\text{max}}^{\frac{1}{2}} \frac{\sqrt{1+a^2+b^2}}{\sqrt{\lambda_{\text{max}} + a^2\lambda_{\text{min}}+ b^2\lambda_{\text{mid}}}} \leq   \lambda_{\text{min}}^{\frac{1}{2}}\lambda_{\text{mid}}^{\frac{1}{2}}\lambda_{\text{max}}^{\frac{1}{2}} \frac{\sqrt{1+a^2+b^2}}{\sqrt{\lambda_{\text{min}} + a^2\lambda_{\text{min}}+ b^2\lambda_{\text{min}}}} = \lambda_{\text{mid}}^{\frac{1}{2}}\lambda_{\text{max}}^{\frac{1}{2}}.
 \en
$\Box$

\subsubsection{Error analysis}
Throughout this paper, we use $C$ to denote a generic constant independent of $A,\omega,h,p$, $u$ and $\hat u$.
The abstract error estimate built in \cite{ref21} also holds in the current situation with the $\big|\big|\big| \cdot
\big|\big|\big|_{\mathcal{F}_h}-$norm.

\begin{lemma}
Let $u$ be the analytical solution of (\ref{3danisohelm})-(\ref{robin}), and let $u_h$ be the approximate solution of (\ref{helmdisvaria}).
Then, we have
 \be \label{helmabstracterror}
 \big|\big|\big| u - u_h \big|\big|\big|_{\mathcal{F}_h}  \leq C
\mathop{\text{inf}}\limits_{ v_h \in { V}_p({\cal T}_h) }
\big|\big|\big| u - v_h \big|\big|\big|_{\mathcal{F}_h^+}.
\en
\end{lemma}

For convenience, we use $\hat V(\hat{\cal T}_{\hat h})$ and $\hat V_p(\hat{\cal T}_{\hat h})$
to denote the images of $V({\cal T}_h)$ and $V_p({\cal T}_h)$ under the coordinate transformation (\ref{tran1}), respectively.
In addition, we define
$\hat{u}_{\hat h}(\hat{\bf x})= {u}_h(S^{-1}\hat{\bf x})$,
and endow the space $\hat V(\hat{\cal T}_{\hat h})$ with the norm
\begin{eqnarray} \label{helmtransfnorm}
&||| \hat{ v}|||_{\hat{\mathcal{F}}_{\hat h}}^2~:=
\omega^{-1}||\beta^{1/2}\llbracket \hat \nabla_{\hat h} \hat v\rrbracket_N||_{0,\hat{\mathcal{F}}_{\hat h}^I}^2
+ \omega||\alpha^{1/2}\llbracket \hat v\rrbracket_N||_{0,\hat{\mathcal{F}}_{\hat h}^I}^2
\cr & + \omega^{-1}||\delta^{1/2}  \hat\nabla_{\hat h} \hat v\cdot \hat{\bf
n}||_{0,\hat{\mathcal{F}}_{\hat h}^B}^2 +
\omega||(1-\delta)^{1/2} \hat v||_{0,\hat{\mathcal{F}}_{\hat h}^B}^2
\end{eqnarray}
and the augmented norm
\begin{eqnarray} \label{helmtransfplusnorm}
& |||\hat v|||_{\hat{\mathcal{F}}_{\hat h}^+}^2~:= |||\hat v|||_{\hat{\mathcal{F}}_{\hat h}}^2 +
\omega||\beta^{-1/2}\{\{ \hat v\}\}||_{0,\hat{\mathcal{F}}_{\hat h}^I}^2 \cr & +
\omega^{-1}||\alpha^{-1/2}\{\{ \hat\nabla_{\hat h} \hat v
\}\}||_{0,\hat{\mathcal{F}}_{\hat h}^I}^2 +
\omega||\delta^{-1/2} \hat v||_{0,\hat{\mathcal{F}}_{\hat h}^B}^2.
\end{eqnarray}

The following Lemma states the transformation stability with respect to a mesh-dependent norm and a mesh-independent norm, respectively.

\begin{lemma} \label{helmimportl} For $u \in V({\cal T}_h)$, we have
\beq \label{helmstade}
\begin{split}
 \big|\big|\big| u \big|\big|\big|_{\mathcal{F}_h}& \leq
~\rho^{\frac{1}{2}} ~\lambda_{\text{mid}}^{\frac{1}{4}}~\lambda_{\text{max}}^{\frac{1}{4}}~ (1 + \lambda_{\text{min}}^{\frac{1}{2}} )~ \big|\big|\big|\hat{u}\big|\big|\big|_{\hat{\mathcal{F}}_{\hat h}}, \\
\big|\big|\big|{u}\big|\big|\big|_{\mathcal{F}_h^+}
&\leq
~\rho^{\frac{1}{2}} ~\lambda_{\text{mid}}^{\frac{1}{4}}~\lambda_{\text{max}}^{\frac{1}{4}}~ (1 + \lambda_{\text{min}}^{\frac{1}{2}} )~\big|\big|\big|\hat{u}\big|\big|\big|_{\hat{\mathcal{F}}_{\hat h}^+},
\end{split}
\eq
and
 \be\label{helmstainde}
 \big|\big| u \big|\big|_{0,\Omega} \leq \bigg(\text{det}(\Lambda^{\frac{1}{2}}) \bigg)^{\frac{1}{2}}
 ~\big|\big|\hat{u}\big|\big|_{0,\hat\Omega}.
\en
\end{lemma}

{\it Proof}.  We divide the proof into three steps.

{\it Step 1:} To estimate
$||\llbracket A \nabla_h u\rrbracket_N||_{0,\mathcal{F}_h^I} $
and $||A \nabla_h u\cdot {\bf n}||_{0,\mathcal{F}_h^B}$.

By the coordinate transformation (\ref{tran1}) and direct calculation, we obtain
\be \label{helmne1}
\nabla u = P^T\Lambda^{-\frac{1}{2}}\hat\nabla \hat u,
\en
which implies that
 \be \label{helmne2}
 A \nabla u = P^T\Lambda^{\frac{1}{2}}\hat\nabla \hat u.
\en
Thus, on the interface $\Gamma_{kj}\in \mathcal{F}_h^I$, we have
 \be \label{helmne3}
 \llbracket A \nabla_h u\rrbracket_N =
 {\bf n}_k \cdot (P^T\Lambda^{\frac{1}{2}}\hat\nabla_{\hat h} {\hat u}_k)
 + {\bf n}_j \cdot (P^T\Lambda^{\frac{1}{2}}\hat\nabla_{\hat h} {\hat u}_j).
\en
It is easy to check that
\be \label{helmne4}
  {\bf n}_k =
|\Lambda^{\frac{1}{2}}P{\bf n}_k|~ P^T\Lambda^{\frac{-T}{2}}\hat{\bf n}_k.
\en
Then
 \be \label{helmne5}
 {\bf n}_k \cdot (P^T\Lambda^{\frac{1}{2}}\hat\nabla {\hat u}_k)
  = |\Lambda^{\frac{1}{2}}P{\bf n}_k|~  \hat{\bf n}_k \cdot \hat\nabla_{\hat h} {\hat u}_k.
 \en
Substituting (\ref{helmne5}) into (\ref{helmne3}), yields
\be \label{helmne7}
\llbracket A \nabla_h u\rrbracket_N
= |\Lambda^{\frac{1}{2}}P{\bf n}_k|~  \llbracket \hat\nabla_{\hat h} {\hat u} \rrbracket_N.
\en
Furthermore, by the scaling argument and \textcolor{red}{ the relation (\ref{arearela}), }
we obtain
\be \label{helmnormaljumpf}
||\llbracket A \nabla_h u\rrbracket_N||_{0,\mathcal{F}_h^I}
\leq ||\Lambda^{\frac{1}{2}}||~ \lambda_{\text{mid}}^{\frac{1}{4}}~\lambda_{\text{max}}^{\frac{1}{4}} ~
||\llbracket \hat\nabla_{\hat h} \hat u \rrbracket_N||_{0,\hat{\mathcal{F}}_{\hat h}^I}
\en
and
 \be \label{helms1_2}
|| A \nabla_h u\cdot {\bf n}||_{0,\mathcal{F}_h^B}
\leq ||\Lambda^{\frac{1}{2}}||~ \lambda_{\text{mid}}^{\frac{1}{4}}~\lambda_{\text{max}}^{\frac{1}{4}} ~
|| \hat\nabla_{\hat h} \hat u ||_{0,\hat{\mathcal{F}}_{\hat h}^B} .
\en

{\it Step 2:} To estimate $\big|\big|\llbracket {u}\rrbracket_N \big|\big|_{0,\mathcal{F}_h^I}$.

By the coordinate transformation (\ref{tran1}), we obtain
\be \label{helmne9}
u_k{\bf n}_k = |\Lambda^{\frac{1}{2}}P{\bf n}_k|~ P^T \Lambda^{\frac{-T}{2}}( {\hat u}_k \hat{\bf n}_k).
 \en
Then, on the interface $\Gamma_{kj}\in \mathcal{F}_h^I$, we have
\be \label{helmnj2}
\llbracket {u}\rrbracket_N = |\Lambda^{\frac{1}{2}}P{\bf n}_k|~
P^T \Lambda^{\frac{-T}{2}} \llbracket \hat{u}\rrbracket_N.
\en
Substituting (\ref{helmnj2}) into $\big|\big|\llbracket {u}\rrbracket_N \big|\big|_{0,\mathcal{F}_h^I}$ and combining with \textcolor{red}{(\ref{arearela}) } and
$|\Lambda^{\frac{1}{2}}P{\bf n}_k|~ ||P^T \Lambda^{\frac{-T}{2}}|| \leq \rho^{\frac{1}{2}}$, yields
\be \label{helmne11}
 \big|\big|\llbracket u \rrbracket_N \big|\big|_{0,\mathcal{F}_h^I} \leq
  \rho^{\frac{1}{2}} ~\lambda_{\text{mid}}^{\frac{1}{4}}~\lambda_{\text{max}}^{\frac{1}{4}}~
  \big|\big|\llbracket \hat{u}\rrbracket_N \big|\big|_{0,\hat{\mathcal{F}}_{\hat h}^I}.
 \en

{\it Step 3:} Build estimates of $||\{\{ u\}\}||_{0,\mathcal{F}_h^I}$,
$||\{\{ A \nabla_h u \}\}||_{0,\mathcal{F}_h^I}$ and $||u||_{0,\mathcal{F}_h^B}$.

By the coordinate transformation (\ref{tran1}) and the scaling argument, we get
\be \label{helmne12}
\big|\big| \{\{ u \}\} \big|\big|_{0,\mathcal{F}_h^I} =
\big|\big| \{\{ \hat u \}\} \big|\big|_{0,\mathcal{F}_h^I}
 ~\leq~ \lambda_{\text{mid}}^{\frac{1}{4}}~\lambda_{\text{max}}^{\frac{1}{4}}~
 \big|\big| \{\{ \hat u \}\} \big|\big|_{0,\hat{\mathcal{F}}_h^I},
\en
\be \label{helmne13}
||\{\{ A\nabla_h u \}\}||_{0,\mathcal{F}_h^I}   \xlongequal{(\ref{helmne2})}
 \big|\big| P^T\Lambda^{\frac{1}{2}} \{\{\hat\nabla_{\hat h} \hat u \}\} \big|\big|_{0,\mathcal{F}_h^I}
\leq ~||\Lambda^{\frac{1}{2}}|| ~ \lambda_{\text{mid}}^{\frac{1}{4}}~\lambda_{\text{max}}^{\frac{1}{4}}~
\big|\big|  \{\{\hat\nabla_{\hat h} \hat u \}\} \big|\big|_{0,\hat{\mathcal{F}}_{\hat h}^I},
\en
and
\begin{eqnarray}\label{helmne14}
& \big|\big|  u \big|\big|_{0,\mathcal{F}_h^B} 
\leq ~\lambda_{\text{mid}}^{\frac{1}{4}}~\lambda_{\text{max}}^{\frac{1}{4}}~
 \big|\big|  \hat u  \big|\big|_{0,\hat{\mathcal{F}}_{\hat h}^B}.
\end{eqnarray}
Combining (\ref{helmnormaljumpf})-(\ref{helms1_2}) with (\ref{helmne11})-(\ref{helmne14}) yields the desired result (\ref{helmstade}).

Finally, by the coordinate transformation (\ref{tran1}) and the scaling argument,
we directly obtain the result (\ref{helmstainde}).

$\Box$

\begin{remark}
As in the proof of the above Lemma, we can obtain the following transformation stability with respect to two mesh-dependent norms,
for $\forall\hat u \in \hat V(\hat{\cal T}_{\hat h})$,
\beq \label{thelmstade}
\begin{split}
 \big|\big|\big| \hat u \big|\big|\big|_{\hat{\mathcal{F}}_{\hat h}} & \leq
~\rho^{\frac{1}{2}} ~\lambda_{\text{min}}^{-\frac{1}{4}}~\lambda_{\text{mid}}^{-\frac{1}{4}}~ (1 + \lambda_{\text{max}}^{-\frac{1}{2}} )~ \big|\big|\big| u \big|\big|\big|_{\mathcal{F}_h}, \\
\big|\big|\big| \hat{u}\big|\big|\big|_{\hat{\mathcal{F}}_{\hat h}^+}  &\leq
~\rho^{\frac{1}{2}} ~\lambda_{\text{min}}^{-\frac{1}{4}}~\lambda_{\text{mid}}^{-\frac{1}{4}}~ (1 + \lambda_{\text{max}}^{-\frac{1}{2}} )~ \big|\big|\big| u \big|\big|\big|_{\mathcal{F}_h^+}.
\end{split}
\eq
\end{remark}

Set $\hat\lambda=\text{min}_{\hat\Omega_k\in{\hat{\cal T}}_{\hat h}}\hat\lambda_k$, where $\hat\lambda_k$ is
the positive parameter depending only on the shape of an element $\hat\Omega_k$
of ${\hat{\cal T}}_{\hat h}$ introduced in \cite[Theorem 3.2]{mhp}.
Let $m$ be a given positive integer satisfying $m\geq 2~(1+2^{1/\hat\lambda})$.

 For a bounded and connected domain \(D\subset \Omega\), let
\(||\cdot||_{s,\omega,D}\) be the \(\omega-\)weighted Sobolev norm
defined by
$$||v||_{s,\omega,D}^2 = \sum_{j=0}^{s}\omega^{2(s-j)}|v|_{j,D}^2. $$

The following Lemma is a direct consequence of Lemma 3.7 in \cite{ref21} and  Corollary 5.5 in \cite{mhp}.

\begin{lemma}\label{helmpwapp}
\textcolor{red}{Let $1< r\leq \frac{m-1}{2}$ with a sufficiently large $m$.}
Assume that the analytical solution $\hat{u}\in H^{r+1}(\hat\Omega)$
satisfies the Helmholtz equation (\ref{helm8}) in isotropic media.
Then there is a function $\hat\xi_{\hat h} \in \hat V_p(\hat{\cal T}_{\hat h})$ such that
 \be \label{helmapprox}
\big|\big|\big| \hat u - \hat\xi_{\hat h} \big|\big|\big|_{\hat{\mathcal{F}}_{\hat h}^+}
\leq ~C_1 ~\omega^{-\frac{1}{2}}~\hat h^{r-\frac{1}{2}} m^{-\hat\lambda(r-\varepsilon)}~||\hat u||_{r+1,\omega,\hat\Omega},
\en
and
\begin{equation} \label{Poincare}
||\hat u  ||_{0,\hat\Omega} \leq
C~(\hat h^{1/2}\omega^{1/2}+\hat h^{-1/2}\omega^{-1/2})~
\big|\big|\big|\hat{u}\big|\big|\big|_{\hat{\mathcal{F}}_{\hat h}},
\end{equation}
where $C_1=C_1(\omega \hat h)$ is independent of $p$ and $\hat u$, but increases as a function of the product $\omega \hat h$ and depends
on the shape of the $\hat\Omega_k\in \hat{\cal T}_{\hat h}$, the index $r$ and the flux parameters.
$\varepsilon=\varepsilon(m)>0$ satisfies $\varepsilon(m)\rightarrow
0$ when $m\rightarrow\infty$.
\end{lemma}

Based on the above discussions, we can easily build the desired error estimates of the approximation $u_h$.
\begin{theorem}  \label{helml2error} Let $u$ and $u_h$ denote the analytical solution of (\ref{3danisohelm})-(\ref{robin}) and the proposed PWDG approximation,
respectively. Suppose that $p$ and $r$ satisfy the conditions in Lemma \ref{helmpwapp}. Assume \textcolor{red}{that $\omega \hat h \leq C$} and $u\in H^{r+1}(\Omega)$.
Then, for sufficiently large $p$, we have
 \be \label{helmfhnormerr}
\big|\big|\big| u - u_{ h} \big|\big|\big|_{{\mathcal{F}}_{ h}^+} \leq
\textcolor{red}{C} ~C_2~\rho^{\frac{3}{4}} ~\omega^{-\frac{1}{2}}~  h^{r-\frac{1}{2}} m^{-\hat\lambda(r-\varepsilon)} ~|| u||_{r+1,\omega,\Omega},
\en
and
\be \label{helm4.estim4}
 ||u-u_h||_{0,\Omega} \leq ~ \textcolor{red}{C}   
 ~C_3~~\rho^{\frac{5}{4}}~\omega^{-1}~ h^{r-1} m^{-\hat\lambda(r-\varepsilon)}~|| u||_{r+1,\omega,\Omega} ,
 \en
 where
 the positive numbers $C_2$ and $C_3$ are defined as
$$ C_2=||A^{\frac{1}{2}}||~ (1 + ||A^{-\frac{1}{2}}||^{-1} )\quad\mbox{and}\quad C_3= (1 + ||A^{-\frac{1}{2}}||^{-1} )~( ||A|| + ||A^{\frac{1}{2}}||). $$
\end{theorem}
{\it Proof}. With the transformations (\ref{tran1}),
we define ${\xi}_{h}({\bf x}) = \hat{\xi}_{\hat h} (S{\bf x})$, where $\hat{\xi}_{\hat h}$ satisfying (\ref{helmapprox})
denotes the plane wave approximation of the scaled acoustic field $\hat{u}$. \textcolor{red}{Under the assumption $\omega \hat h \leq C$, we have $C_1(\omega \hat h)\leq C$.}
Thus, using (\ref{helmabstracterror}) and (\ref{helmstade}), together with (\ref{helmapprox}), leads to
\begin{eqnarray}  \label{inteerr1}
 & \big|\big|\big| u - u_h \big|\big|\big|_{\mathcal{F}_h}
 \leq
  ~C~\big|\big|\big| u - {\xi}_{ h} \big|\big|\big|_{{\mathcal{F}}_{ h}^+}
  \leq
  ~C~\rho^{\frac{1}{2}}  ~\lambda_{\text{mid}}^{\frac{1}{4}}~\lambda_{\text{max}}^{\frac{1}{4}}~ (1 + \lambda_{\text{min}}^{\frac{1}{2}} )~ \big|\big|\big| \hat{u} - \hat{\xi}_{\hat h} \big|\big|\big|_{\hat{\mathcal{F}}_{\hat h}^+}
  \cr & \leq
   ~C~\rho^{\frac{1}{2}} ~\lambda_{\text{mid}}^{\frac{1}{4}}~\lambda_{\text{max}}^{\frac{1}{4}}~ (1 + \lambda_{\text{min}}^{\frac{1}{2}} )~ \omega^{-\frac{1}{2}}~ \hat h^{r-\frac{1}{2}} m^{-\hat\lambda(r-\varepsilon)} ||\hat u||_{r+1,\omega,\hat\Omega}.
\end{eqnarray}
    Then, by the scaling argument and (\ref{geometric100}), we further obtain
\beq
\big|\big|\big| u - u_h \big|\big|\big|_{\mathcal{F}_h}
 &\leq& ~C~\rho^{\frac{1}{2}} ~\omega^{-\frac{1}{2}}~h^{r-\frac{1}{2}}~ m^{-\hat\lambda(r-\varepsilon)}
 ~\lambda_{\text{mid}}^{\frac{1}{4}}~\lambda_{\text{max}}^{\frac{1}{4}}~ (1 + \lambda_{\text{min}}^{\frac{1}{2}} )~ ||\Lambda^{\frac{1}{2}}||^{\frac{3}{2}}~ (det(\Lambda^{-\frac{1}{2}}))^{\frac{1}{2}} || u||_{r+1,\omega,\Omega}\cr
 & \leq& ~C~||\Lambda^{\frac{1}{2}}||~(1 + ||\Lambda^{-\frac{1}{2}}||^{-1} )~\rho^{\frac{3}{4}}~ \omega^{-\frac{1}{2}}~h^{r-\frac{1}{2}}~ m^{-\hat\lambda(r-\varepsilon)}
 || u||_{r+1,\omega,\Omega}.
 \label{fhnorm100}
\eq
Combining (\ref{helmstainde}), (\ref{Poincare}), (\ref{geometric100}), (\ref{thelmstade}) and (\ref{fhnorm100}),
yields
 \beq \nonumber 
\begin{split}
\big|\big| u-u_h \big|\big|_{0,\Omega}
& \overset{(\ref{helmstainde})} {\leq}  \bigg(\text{det}(\Lambda^{\frac{1}{2}}) \bigg)^{\frac{1}{2}}~
\big|\big| \hat{u}- \hat u_{\hat h} \big|\big|_{0,\hat\Omega}&
\\ &  \overset{(\ref{Poincare})} {\leq}
C~ \bigg(\text{det}(\Lambda^{\frac{1}{2}}) \bigg)^{\frac{1}{2}} ~ ( \omega^{-\frac{1}{2}}\hat h^{-\frac{1}{2}} + \omega^{\frac{1}{2}}\hat h^{\frac{1}{2}} ) ~
\big|\big|\big|\hat{u} - \hat u_{\hat h}  \big|\big|\big|_{\hat{\mathcal{F}}_{\hat h}}  &
\\ & \overset{(\ref{geometric100})}{ \leq} C ~ \bigg(\text{det}(\Lambda^{\frac{1}{2}}) \bigg)^{\frac{1}{2}}  ~ \omega^{-\frac{1}{2}} h^{-\frac{1}{2}} ~||\Lambda^{\frac{1}{2}}||^{\frac{1}{2}} ~ \big|\big|\big|\hat{u} - \hat u_{\hat h}  \big|\big|\big|_{\hat{\mathcal{F}}_{\hat h}}  &
\\ & \overset{(\ref{thelmstade})} {\leq} C ~ \bigg(\text{det}(\Lambda^{\frac{1}{2}}) \bigg)^{\frac{1}{2}} ~ \omega^{-\frac{1}{2}}  ~h^{-\frac{1}{2}} ~ ||\Lambda^{\frac{1}{2}}||^{\frac{1}{2}} ~\rho^{\frac{1}{2}} ~\lambda_{\text{min}}^{-\frac{1}{4}}~\lambda_{\text{mid}}^{-\frac{1}{4}}~ (1 + \lambda_{\text{max}}^{-\frac{1}{2}} )~\big|\big|\big| u - u_h  \big|\big|\big|_{\mathcal{F}_{ h}}  &
\\ &  \overset{(\ref{fhnorm100})} {\leq} ~ C ~C_3~\rho^{\frac{5}{4}} ~ \omega^{-1} ~ h^{r-1} ~m^{-\hat\lambda(r-\varepsilon)}
 || u||_{r+1,\omega,\Omega},  & 
\end{split}
\eq
 This completes the proof.

 $\Box$

\begin{remark} \label{helml2comp}
We emphasize that \textcolor{red}{the shape regularity assumption} on the triangulation $\hat{\mathcal T}_{\hat{h}}$ (instead of ${\cal T}_h$) is important in the
derivation of the error estimate (\ref{helm4.estim4}). In fact, if we directly assume that the triangulation ${\cal T}_h$ is shape regular, we have different relations
on the two mesh sizes:
\be
||\Lambda^{\frac{1}{2}}||^{-1}h \leq \hat h \leq ||\Lambda^{-\frac{1}{2}}||~h.
\en
In this situation, we can only build a weaker $L^2$ error estimate than (\ref{helm4.estim4})
\be
 ||u-u_h||_{0,\Omega} \leq ~C~ C_3 ~\rho^{\frac{r}{2}+1}~ \omega^{-1} ~ h^{r-1} m^{-\theta\lambda(r-\varepsilon)}
 || u||_{r+1,\omega,\Omega}.\quad\quad(r>1)
\en
\end{remark}

\begin{remark} A natural idea is to apply the standard PWDG method to the isotropic Helmholtz equation (\ref{helm8}) derived by the transformation $S$
and then use the image of the resulting approximation under the inverse transformation $S^{-1}$ as the desired approximation of $u$, but this idea
may be disappointing. Let $\hat {u}(\hat{\bf x})$ and $\tilde{u}_{\hat h}(\hat{\bf x})$ denote the analytic solution of the
equation (\ref{helm8}) and its PWDG approximation, respectively, and let $\tilde{u}_h({\bf x})$ denote the image of $\tilde{u}_{\hat h}(\hat{\bf x})$ under the
inverse transformation $S^{-1}$. Although the following abstract error estimate is still valid
\be \label{helmHATabstracterror}
\big|\big|\big| \hat{u} - \tilde{u}_{\hat h} \big|\big|\big|_{\hat{\mathcal{F}}_{\hat h}}
\leq  \mathop{\text{inf}}\limits_{\hat{ \xi}_{\hat h}\in \hat{ V}_p(\hat{\cal T}_{\hat h})}
\big|\big|\big| \hat{u} - \hat{ \xi}_{\hat h} \big|\big|\big|_{\hat{\mathcal{F}}_{\hat h}^+},
\en
the approximation $\tilde{u}_h$ does not satisfy the abstract error estimate (\ref{helmabstracterror}).
Then it seems impossible to build a desired error estimate of the approximation $\tilde{u}_h$.
In fact, by (\ref{helmne7}) and (\ref{helmnj2}), we have
$$\llbracket A \nabla_h v\rrbracket_N= |\Lambda^{\frac{1}{2}}P{\bf n}_k|~  \llbracket \hat\nabla_{\hat h} {\hat v} \rrbracket_N
\quad\mbox{and}\quad \llbracket {v}\rrbracket_N = |\Lambda^{\frac{1}{2}}P{\bf n}_k|~
P^T \Lambda^{\frac{-T}{2}} \llbracket \hat{v}\rrbracket_N\quad \mbox{on}~\Gamma_{kj}\in \mathcal{F}_h^I. $$
\textcolor{red}{The different factors in the right hands of these relations} tell us that the proposed PWDG approximation $u_h$ is indeed different from
the image $\tilde{u}_h({\bf x})$ of $\tilde{u}_{\hat h}(\hat{\bf x})$ under the inver transformation $S^{-1}$.
\end{remark}

\section{Plane wave method for three-dimensional anisotropic time-harmonic Maxwell's equations}
In this section we extend the method proposed in the last section to time-harmonic Maxwell's equations in
three-dimensional anisotropic media. As we will see, the current situation is more complex than the case of Helmholtz equation.

\subsection{The model and its variational formula}
We consider
three-dimensional time-harmonic Maxwell equations written as a first-order system of
equations:
\begin{equation} \label{eq1}
\left\{ \begin{aligned}
     &   \nabla\times{\bf E} - i\omega\mu{\bf H}={\bf 0} \\
      & \nabla\times{\bf H} +i \omega\varepsilon{\bf E}={\bf 0} \\
      &  \nabla \cdot (\varepsilon{\bf E}) = 0 \\
      &   \nabla \cdot (\mu{\bf H}) = 0
                          \end{aligned} \right.  \quad {\text in} \quad \Omega
                          \end{equation}
with the lowest-order absorbing boundary condition
 \begin{equation} \label{eq2}
{\bf H}\times{\bf n}-\vartheta({\bf n}\times{\bf E})\times{\bf n}=
{\bf g}/i\omega \quad \text{on} \quad\gamma=\partial\Omega.
 \end{equation}
Here ${\bf E}=(E_x,E_y,E_z)^T$, ${\bf H}=(H_x,H_y,H_z)^T$; \(\omega>0\) is the
temporal frequency of the field; $\vartheta \neq 0$ is assumed to be
constant; \( {\bf g}\in L^2(\partial\Omega)\).
The permittivity \(\varepsilon\) and the permeability
\(\mu\) are assumed to be of the form
\begin{eqnarray} \label{standa1}
& \varepsilon = \varepsilon_r A, \cr & \mu = \mu_r A.
\end{eqnarray}
where $\varepsilon_r, \mu_r$ are constant, and $A$ is assumed to be real strictly positive definite matrix.

We also denote by \( {\cal T}_h\) the partition of the domain $\Omega$, which is also a bounded polyhedron.
Define the broken Sobolev space
\be \label{sobo1}
{\bf H}^r(\text{curl};{\cal T}_h)=\{ {\bf w}\in L^2(\Omega)^3:{\bf
w}|_{\Omega_k}\in {\bf H}^r(\text{curl};\Omega_k) \quad \quad \forall
\Omega_k\in {\cal T}_h \}.
\en
Let ${\bf V}({\cal T}_h)$ be the piecewise Trefftz space defined on ${\cal T}_h$ by
\be
\begin{split}
&  {\bf V}({\cal T}_h) = \bigg\{ {\bf w}\in L^2(\Omega)^3: \exists s>0
~s.t.~ {\bf w} \in {\bf H}^{1/2+s}(\text{curl};{\cal T}_h), \\
& \quad\quad\quad \text{and}~\nabla\times(\mu^{-1}\nabla\times{\bf w} )
-\omega^2\varepsilon {\bf w} ={\bf 0} \text{~in ~ each } ~ \Omega_k\in
{\cal T}_h \bigg\}.  \\
\end{split}
\en

Let $\alpha,\beta,\delta$ be strictly positive constants, with $0<\delta\leq 1/2$.
Define the sesquilinear form $\mathcal{A}_h(\cdot, \cdot)$ by
\beqx
\mathcal{A}_h({\bf E}, {\bm \xi}) &=&  -\int_{\mathcal{F}_h^{\text{I}}} \{\{{\bf E}\}\}\cdot
\overline{\llbracket \mu^{-1} \nabla_h\times {\bm
\xi}\rrbracket_T}~dS -i\omega^{-1}\int_{\mathcal{F}_h^{\text{I}}}
\beta \llbracket \mu^{-1} \nabla_h\times {\bf E}\rrbracket_T \cdot
\overline{\llbracket \mu^{-1} \nabla_h\times
{\bm\xi}\rrbracket_T}~dS
\\
&-& \int_{\mathcal{F}_h^{\text{I}}} \{\{ \mu^{-1}
\nabla_h\times {\bf E}\}\}\cdot \overline{\llbracket
{\bm\xi}\rrbracket_T}~dS -
i\omega\int_{\mathcal{F}_h^{\text{I}}}\alpha \llbracket {\bf
E}\rrbracket_T \cdot \overline{\llbracket {\bm\xi}\rrbracket_T}~dS
\\
&+&\int_{\mathcal{F}_h^{\text{B}}}(1-\delta)({\bf n}\times{\bf
E})\cdot\overline{(\mu^{-1} \nabla_h\times {\bm\xi})}~dS
- \int_{\mathcal{F}_h^{\text{B}}} \delta(\mu^{-1} \nabla_h\times {\bf
E})\cdot\overline{({\bf n}\times{\bm\xi})}~dS
\\
&-&
i\omega^{-1}\int_{\mathcal{F}_h^{\text{B}}}\delta\vartheta^{-1}[{\bf
n}\times(\mu^{-1} \nabla_h\times {\bf E})] \cdot \overline{[{\bf
n}\times(\mu^{-1} \nabla_h\times {\bm\xi})]}~dS
\\
&-&
i\omega\int_{\mathcal{F}_h^{\text{B}}} (1-\delta)\vartheta({\bf
n}\times {\bf E}) \cdot \overline{({\bf n}\times{\bm\xi})}~dS, ~\forall {\bm \xi}\in {\bf V}({\cal T}_h)
\eqx
and the functional $\ell_h(\cdot,\cdot)$ by
\be \label{nonhomovari3}
  \ell_h({\bf g},{\bm \xi})=
-i\omega^{-1}\int_{\mathcal{F}_h^{\text{B}}}\delta\vartheta^{-1}({\bf
n}\times {\bf g})\cdot \overline{\mu^{-1} \nabla_h\times {\bm
\xi}}~dS + \int_{\mathcal{F}_h^{\text{B}}} (1-\delta)({\bf n}\times
{\bf g})\cdot \overline{({\bf n}\times{\bm\xi})}~dS,~\forall {\bm \xi}\in {\bf V}({\cal T}_h)
 \en

\textcolor{red}{Then, for a given ${\bf g}$, the variational problem associated with (\ref{eq1})-(\ref{eq2}) can be expressed as follows (see \cite[Section 3]{pwdg}). Find ${\bf
E}\in {\bf V}({\cal T}_h)$ such that,}
\be \label{nonhomovari2}
\mathcal{A}_h({\bf E}, {\bm \xi}) =\ell_h({\bf g},{\bm \xi}),
\quad \forall {\bm \xi} \in {\bf V}({\cal T}_h).
\en

\subsection{Plane wave discretization for the 3D anisotropic Maxwell equations}
The proposed plane wave method for (\ref{eq1}) depends on two transformations.

\subsubsection{A scaled transformation and a coordinate transformation}
\textcolor{red}{
Since $A$ is positive definite matrix, there exists an orthogonal matrix $P$ and a diagonal positive definite matrix $\Lambda = \text{diag} ( \lambda_{\text{min}},\lambda_{\text{mid}},\lambda_{\text{max}} )$ such that
$A= P^T \Lambda P$, where $\lambda_{\text{min}} \leq \lambda_{\text{mid}} \leq \lambda_{\text{max}}$ are constant. Without loss of generality, we also assume that $\text{det}(P)=1$.
Furthermore, we set
\be  \label{mdefi}
m_{\text{max}} = \sqrt{\lambda_{\text{mid}}\lambda_{\text{max}}}, \quad m_{\text{mid}} = \sqrt{\lambda_{\text{max}}\lambda_{\text{min}}}\quad\mbox{and}\quad m_{\text{min}} =
\sqrt{\lambda_{\text{min}}\lambda_{\text{mid}}}.
\en
}

Define the scaled fields $\tilde{\bf E}$ and $\tilde{\bf H}$ as
\begin{eqnarray}
\label{standa2}
& (E_x,E_y,E_z)^T = G~(\tilde{E}_x,\tilde{E}_y,\tilde{E}_z)^T,
\cr &
  (H_x,H_y,H_z)^T = G~(\tilde{H}_x,\tilde{H}_y,\tilde{H}_z)^T.
\end{eqnarray}
Here $G = P^T\Lambda^{-\frac{1}{2}}$.
\textcolor{red}{Then, by direct calculation on (\ref{standa2}), we deduce that}
\begin{eqnarray} \label{transrela1}
& \nabla \times {\bf E} = P^T\Lambda^{\frac{1}{2}}
( -\frac{q_3\cdot \nabla \tilde{E}_y}{m_{\text{min}}} + \frac{q_2\cdot \nabla \tilde{E}_z}{m_{\text{mid}}},
 \frac{q_3\cdot \nabla \tilde{E}_x}{m_{\text{min}}} - \frac{q_1\cdot \nabla \tilde{E}_z}{m_{\text{max}}},
 -\frac{q_2\cdot \nabla \tilde{E}_x}{m_{\text{mid}}} + \frac{q_1\cdot \nabla \tilde{E}_y}{m_{\text{max}}})^T,
\cr &
\nabla \times {\bf H} = P^T\Lambda^{\frac{1}{2}}
( -\frac{q_3\cdot \nabla \tilde{H}_y}{m_{\text{min}}} + \frac{q_2\cdot \nabla \tilde{H}_z}{m_{\text{mid}}},
 \frac{q_3\cdot \nabla \tilde{H}_x}{m_{\text{min}}} - \frac{q_1\cdot \nabla \tilde{H}_z}{m_{\text{max}}},
 -\frac{q_2\cdot \nabla \tilde{H}_x}{m_{\text{mid}}} + \frac{q_1\cdot \nabla \tilde{H}_y}{m_{\text{max}}})^T,
 \cr &
 \mu{\bf H} = \mu_r P^T \Lambda P P^T\Lambda^{-\frac{1}{2}}\tilde{\bf H} = \mu_r P^T\Lambda^{\frac{1}{2}}\tilde{\bf H},
 \cr &
 \varepsilon {\bf E} = \varepsilon_r P^T \Lambda P P^T\Lambda^{-\frac{1}{2}}\tilde{\bf E} = \varepsilon_r P^T\Lambda^{\frac{1}{2}}\tilde{\bf E}.
\end{eqnarray}
\textcolor{red}{Here we need to first use (\ref{standa1}) for the derivation of the later two equalities}.

Set $M = \text{diag}(m_{\text{max}},m_{\text{mid}},m_{\text{min}})$, and define the coordinate transformation
\be \label{trans}
 \hat{\bf x}=(\hat{x} ~ \hat{y} ~ \hat{z})^T = ~M~P~(x~y~z)^T\xlongequal{\Delta} S~{\bf x}, \quad S=M~P. 
\en
With the inverse transformation $S^{-1}$, we define the scaled  electric and magnetic fields
\be \label{codimax}
(\hat{\bf E}(\hat{\bf x}), \hat{\bf H}(\hat{\bf x}))^T=(\tilde{\bf E}(S^{-1}\hat{\bf x}), \tilde{\bf H}(S^{-1}\hat{\bf x}))^T=(\tilde{\bf E}(x), \tilde{\bf H}(x))^T.
\en
\textcolor{red}{By direct manipulation on (\ref{codimax}), we obtain the following system }
\begin{eqnarray} \label{transrela16}
& \hat\nabla\times\hat{\bf E}=( -\frac{q_3\cdot \nabla \tilde{E}_y}{m_{\text{min}}} + \frac{q_2\cdot \nabla \tilde{E}_z}{m_{\text{mid}}},
 \frac{q_3\cdot \nabla \tilde{E}_x}{m_{\text{min}}} - \frac{q_1\cdot \nabla \tilde{E}_z}{m_{\text{max}}},
 -\frac{q_2\cdot \nabla \tilde{E}_x}{m_{\text{mid}}} + \frac{q_1\cdot \nabla \tilde{E}_y}{m_{\text{max}}})^T,
 \cr &
\hat\nabla\times\hat{\bf H}=(-\frac{q_3\cdot \nabla \tilde{H}_y}{m_{\text{min}}} + \frac{q_2\cdot \nabla \tilde{H}_z}{m_{\text{mid}}},
 \frac{q_3\cdot \nabla \tilde{H}_x}{m_{\text{min}}} - \frac{q_1\cdot \nabla \tilde{H}_z}{m_{\text{max}}},
 -\frac{q_2\cdot \nabla \tilde{H}_x}{m_{\text{mid}}} + \frac{q_1\cdot \nabla \tilde{H}_y}{m_{\text{max}}})^T,
 \cr &
 M \Lambda^{\frac{1}{2}} \hat\nabla\cdot \hat{\bf H}=\nabla\cdot( P^T\Lambda^{\frac{1}{2}}\tilde{\bf H}),
 \cr &
 M \Lambda^{\frac{1}{2}} \hat\nabla\cdot \hat{\bf E}=\nabla\cdot( P^T\Lambda^{\frac{1}{2}}\tilde{\bf E}).
\end{eqnarray}
\textcolor{red}{Thus, by (\ref{eq1}), (\ref{transrela1}) and (\ref{transrela16}),  the scaled  electric and magnetic fields $(\hat{\bf E}(\hat{\bf x}), \hat{\bf H}(\hat{\bf x}))$
satisfy the transformed isotropic Maxwell equations:}
\begin{equation} \label{ntranseq1}
\left\{ \begin{aligned}
     &   \hat\nabla\times\hat{\bf E} - i\omega\mu_r\hat{\bf H}={\bf 0} \\
      & \hat\nabla\times\hat{\bf H} +i \omega\varepsilon_r\hat{\bf E}={\bf 0} \\
      &  \hat\nabla \cdot (\varepsilon_r \hat{\bf E}) = 0 \\
      &   \hat\nabla \cdot (\mu_r \hat{\bf H}) = 0
                          \end{aligned} \right.  \quad {\text in} \quad \hat\Omega.
                          \end{equation}

Conversely, if the scaled  electric and magnetic fields $(\hat{\bf E}(\hat{\bf x}), \hat{\bf H}(\hat{\bf x}))$
satisfy the transformed isotropic Maxwell equations (\ref{ntranseq1}),
the physical electromagnetic fields $\big({\bf E}({\bf x}),{\bf H}({\bf x})\big)$
\begin{equation} \label{stand7}
\left\{ \begin{aligned}
     &  {\bf E}({\bf x})=  G \hat{\bf E}(\hat{\bf x}) = G \hat{\bf E} (S{\bf x})  \\
      & {\bf H}({\bf x}) =  G \hat{\bf H}(\hat{\bf x}) = G \hat{\bf H} (S{\bf x})
                          \end{aligned} \right.
                          \end{equation}
satisfy the original anisotropic Maxwell equations (\ref{eq1}).

As in Section 2, let $\hat{\Omega}$ and $\hat{\Omega}_k$ denote
the images of $\Omega$ and $\Omega_k$ under the coordinate transformation (\ref{trans}), respectively.
In addition, let $\hat{{\cal T}}_{\hat{h}}$ denote the
partition comprised of the elements
$\{\hat{\Omega}_k\}$, where $\hat{h}$ is the mesh size of the
partition $\hat{{\cal T}}_{\hat{h}}$. Set
$\hat{\mathcal{F}}_{\hat h} = \bigcup_k\partial\hat\Omega_k$,
$\hat{\mathcal{F}}_{\hat h}^{\text{B}}=\partial\hat{\Omega}$ and
$\hat{\mathcal{F}}_{\hat h}^{\text{I}}= \hat{\mathcal{F}}_{\hat h}
\backslash \hat{\mathcal{F}}_{\hat h}^{\text{B}}$.

In the next subsubsection we present a plane wave discretization method for the considered anisotropic Maxwell equations
by using the scaling matrix (\ref{standa2}) and the coordinate transformation (\ref{trans}).

\subsubsection{Anisotropic plane wave basis function spaces}

The discretization is based on a finite-dimensional subspace
 \( {\bf V}_p({\cal T}_h)  \subset {\bf V}({\cal T}_h)\). We first define
plane wave basis function space $\hat{\bf V}_p(\hat{\cal T}_{\hat h})$ satisfying the isotropic Maxwell equations (\ref{ntranseq1}).

By choosing \(p\) unit propagation directions \(
{\bf d}_{l}~(l=1,\cdots,p) \), which can be determined by the codes in \cite{refsite}, we can define plane wave
functions \( \hat{\bf E}_{l}\):
\begin{equation}
\hat{\bf E}_{l}=~{\bf F}_{l}~\text{exp}(\text{i}~\kappa~{\bf
d}_{l}\cdot \hat{\bf x})\quad \text{and} \quad \hat{\bf
E}_{l+p}=~{\bf G}_{l}~\text{exp}(\text{i}~\kappa~{\bf d}_{l}\cdot
\hat{\bf x})~~~(l=1,\cdots,p), \label{maxd1}
\end{equation}
where $\kappa = \omega\sqrt{\mu_r\varepsilon_r}$, ${\bf F}_{l}$ and
${\bf G}_{l}$ are polarization vectors satisfying ${\bf F}_{l}\cdot
{\bf d}_{l}=0$ and ${\bf G}_{l}={\bf F}_{l}\times{\bf d}_{l}~
(l=1,\cdots,p)$.

Let \( \hat{\cal Q}_{2p}\) denote the space spanned by the \(2p\)
plane wave functions \( \hat{\bf E}_{l}\) (\(l=1,\cdots,2p\)),
and define the isotropic plane wave space
\begin{equation}
\hat{\bf V}_p(\hat{\cal T}_{\hat h})=\bigg\{\hat{\bf v}\in
L^2(\hat\Omega):~\hat{\bf v}|_{\hat K}\in \hat{\cal
Q}_{2p}~~\mbox{for~~any}~~\hat K\in \hat{\cal T}_{\hat h} \bigg\}.
\label{hatmaxd2}
\end{equation}

 By (\ref{stand7}), we can define the anisotropic plane wave basis
functions satisfying the original equations (\ref{eq1}):
\begin{equation}\label{dispw}
{\bf E}_{l}=~G~{\bf F}_{l}~\text{exp}(\text{i}~\kappa~{\bf
d}_{l}\cdot S{\bf x})\quad \text{and} \quad {\bf E}_{l+p}=~G~{\bf
G}_{l}~\text{exp}(\text{i}~\kappa~{\bf d}_{l}\cdot S{\bf
x})~~~(l=1,\cdots,p).
\end{equation}

Let \( {\cal Q}_{2p}\) denote the space spanned by the \(2p\)
plane wave basis functions \( {\bf E}_{l}\) (\(l=1,\cdots,2p\)), and define
the anisotropic plane wave space
\begin{equation}
{\bf V}_p({\cal T}_h)=\bigg\{{\bf v}\in L^2(\Omega):~{\bf v}|_K\in
{\cal Q}_{2p}~~\mbox{for~~any}~~K\in {\cal T}_h \bigg\}.
\label{maxd2}
\end{equation}
It is clear that the above space has \(N\times 2p\) basis
functions, which are defined by
\begin{eqnarray}
\bm\phi^k_{l}({\bf x})=\left\{\begin{array}{ll} {\bf E}_{l}({\bf
x}),~~{\bf
x}\in\Omega_k,\\
 0,~~{\bf
x}\in\Omega_j~~\mbox{satisfying}~~j\neq k~~
\end{array}\right.~~(k=1,\cdots,N;~l=1,\cdots,2p).
\label{maxd3}
\end{eqnarray}

Furthermore, for a given ${\bf g}$, we obtain the discretized version of the continuous variational problem (\ref{nonhomovari2}):
 Find ${\bf E}_h \in {\bf V}_p({\cal T}_h)$ such that,
\be \label{maxwelldis}
\mathcal{A}_h({\bf E}_h, {\bm \xi}_h) =\ell_h({\bf g},{\bm \xi}_h),
\quad \forall {\bm \xi}_h \in {\bf V}_p({\cal T}_h).
\en

For convenience, we use $\hat {\bf V}(\hat{\cal T}_{\hat h})$ to denote the image space of the space ${\bf V}({\cal T}_h)$ under the scaling transformation (\ref{standa2}) and the coordinate transformation (\ref{trans}),
and set $\hat{\bf E}_{\hat h}(\hat{\bf x})= G^{-1}{\bf E}_h(S^{-1}\hat{\bf x})$.

\subsection{Error estimates of the approximate solutions}

For the simplicity of notation, we also use $\rho$ to denote the condition number $\text{cond}(A)$ for the positive definite matrix $A$ defined in (\ref{standa1}).
It is easy to see that $\text{cond}(A)= \text{cond}(\Lambda)$ and $\text{cond}(S)=\text{cond}(M) = \rho^{\frac{1}{2}}$.

We adopt similar steps described in section 2 to partition $\Omega$. The minor difference is that,
under the coordinate transformation (\ref{trans}), we have
\be \label{relahhath}
c_0 ||M^{-1}||^{-1}h \leq \hat{h} \leq C_0||M^{-1}||^{-1}h.
\en

We endow ${\bf V}({\cal T}_h)$ with the mesh-skeleton norm,
\be \label{treffznorm}  \begin{split}
&\big|\big|\big|{\bf w}\big|\big|\big|_{\mathcal{F}_h}^2 =
\omega^{-1}\big|\big|\beta^{1/2}\llbracket \mu^{-1}\nabla_h\times
{\bf w}\rrbracket_T\big|\big|_{0,\mathcal{F}_h^I}^2 +
\omega\big|\big|\alpha^{1/2}\llbracket {\bf w}\rrbracket_T
\big|\big|_{0,\mathcal{F}_h^I}^2
\\
&\quad\quad\quad+\omega^{-1}\big|\big|\delta^{1/2}\vartheta^{-1/2} {\bf
n}\times(\mu^{-1}\nabla_h\times {\bf
w})\big|\big|_{0,\mathcal{F}_h^B}^2 +
\omega\big|\big|(1-\delta)^{1/2}\vartheta^{1/2} {\bf n}\times {\bf
w}\big|\big|_{0,\mathcal{F}_h^B}^2\\
\end{split}\en
and the following augmented norm 
\be \label{treffzaugnorm}  \begin{split}
&\big|\big|\big|{\bf w}\big|\big|\big|_{\mathcal{F}_h^+}^2 =
\big|\big|\big|{\bf w}\big|\big|\big|_{\mathcal{F}_h}^2 + \omega
\big|\big| \beta^{-1/2}\{\{ {\bf w} \}\}
\big|\big|_{0,\mathcal{F}_h^I}^2
\\
 &\quad\quad\quad + \omega^{-1} \big|\big|
\alpha^{-1/2}\{\{ \mu^{-1}\nabla_h\times {\bf w} \}\}
\big|\big|_{0,\mathcal{F}_h^I}^2 + \omega \big|\big|
\delta^{-1/2}\vartheta^{1/2} ( {\bf n}\times{\bf w} )
\big|\big|_{0,\mathcal{F}_h^B}^2.\\
\end{split}\en

Similarly, we endow $\hat{\bf V}(\hat{\cal T}_{\hat h})$ with the mesh-skeleton norm, which were introduced in \cite{pwdg},
\be \label{treffztransnorm}
\begin{split}
&\big|\big|\big|\hat{\bf w}\big|\big|\big|_{\hat{\mathcal{F}}_{\hat h}}^2 =
\omega^{-1}\big|\big|\beta^{1/2}\llbracket \mu_r^{-1}\hat \nabla_{\hat h}\times
\hat{\bf w}\rrbracket_T\big|\big|_{0,\hat{\mathcal{F}}_{\hat h}^I}^2 +
\omega\big|\big|\alpha^{1/2}\llbracket \hat{\bf w}\rrbracket_T
\big|\big|_{0,\hat{\mathcal{F}}_{\hat h}^I}^2
\\
&\quad\quad\quad +\omega^{-1}\big|\big|\delta^{1/2}\vartheta^{-1/2}
 \hat{\bf n}\times(\mu_r^{-1} \hat\nabla_{\hat h}\times \hat{\bf
w})\big|\big|_{0,\hat{\mathcal{F}}_{\hat h}^B}^2 +
\omega\big|\big|(1-\delta)^{1/2}\vartheta^{1/2} \hat{\bf n}\times \hat{\bf
w}\big|\big|_{0,\hat{\mathcal{F}}_{\hat h}^B}^2\\
\end{split}
\en
and the following augmented norm
\be \label{treffztransaugnorm}  \begin{split}
&\big|\big|\big| \hat{\bf w}\big|\big|\big|_{\hat{\mathcal{F}}_{\hat h}^+}^2 =
\big|\big|\big| \hat{\bf w}\big|\big|\big|_{\hat{\mathcal{F}}_{\hat h}}^2 + \omega
\big|\big| \beta^{-1/2}\{\{ \hat{\bf w} \}\}
\big|\big|_{0,\hat{\mathcal{F}}_{\hat h}^I}^2
\\
 &\quad\quad\quad + \omega^{-1} \big|\big|
\alpha^{-1/2}\{\{ \mu_r^{-1} \hat\nabla_{\hat h}\times \hat{\bf w} \}\}
\big|\big|_{0,\hat{\mathcal{F}}_{\hat h}^I}^2 + \omega \big|\big|
\delta^{-1/2}\vartheta^{1/2} ( \hat{\bf n}\times \hat{\bf w} )
\big|\big|_{0,\hat{\mathcal{F}}_{\hat h}^B}^2.\\
\end{split}\en

As in \cite{pwdg}, we can show the following
existence, uniqueness and continuity results of solution of the above variational problem.
\begin{lemma}\label{existetc}
There exists a unique solution ${\bf E}_h$ of (\ref{nonhomovari2}).
Moreover, we have \begin{eqnarray} \label{uniquebound} &
-\text{Im}[\mathcal{A}_h({\bf w},{\bf w})] = \big|\big|\big|{\bf
w}\big|\big|\big|^2_{\mathcal{F}_h},
 \cr & \text{and}
~~\big|\mathcal{A}_h({\bf w},{\bm \xi})\big| \leq 2
\big|\big|\big|{\bf w}\big|\big|\big|_{\mathcal{F}_h^+}~
\big|\big|\big|{\bm \xi}\big|\big|\big|_{\mathcal{F}_h},
 ~~\forall ~{\bf w},{\bm \xi}\in {\bf V}({\cal T}_h).
\end{eqnarray}
\end{lemma}

The abstract error estimate built in \cite{pwdg} also holds in the current situation with the \\
$\big|\big|\big| \cdot \big|\big|\big|_{\mathcal{F}_h}-$norm.
\begin{lemma}
Let ${\bf E}$ be the analytical solution of (\ref{eq1})-(\ref{eq2}) and let ${\bf E}_h$ be the approximate solution of (\ref{maxwelldis}). Then, there exists a constant $C$ independent of $\omega,h,p$ and $A$ such that
\be \label{abstracterror}
\big|\big|\big| {\bf
E} - {\bf E}_h \big|\big|\big|_{\mathcal{F}_h}  \leq 3
\mathop{\text{inf}}\limits_{{\bm \xi}_{h}\in {\bf V}_p({\cal T}_h)
} \big|\big|\big| {\bf E} - {\bm \xi}_h
\big|\big|\big|_{\mathcal{F}_h^+}.
\en
\end{lemma}

Let $m$ be a positive integer satisfying the condition described in Subsection 2.3.2, assume that $1<r\leq {m-1\over 2}$ and choose $p=(m+1)^2$.
The following Lemma is a direct consequence of Corollary 5.5 in \cite{pwdg}.

\begin{lemma} \label{approesti}
Assume that the analytical solution $\hat{\bf E}\in
{\bf H}^{r+1}(\text{curl};\hat\Omega)~$
 satisfies the Maxwell equations (\ref{ntranseq1}) in isotropic media.
Choose a set of $p$ plane wave
propagation directions $\{ {\bf d}_l \}_{1\leq l\leq p}$ with the
corresponding set of polarization directions ${\bf  F}_l, {\bf G}_l$
defined by (\ref{maxd1}) in a suitable manner.
Then, 
there is a function $\hat{\bm \xi}_{\hat h}\in \hat{\bf V}_p(\hat{\cal T}_{\hat h})$ such that, for large $p$,
 \be \label{approesti2}
 \big|\big|\big| \hat{\bf E}-\hat{\bm \xi}_{\hat h}
\big|\big|\big|_{\hat{\mathcal{F}}_{\hat h}^+} \leq
C_1 ~\omega^{-5/2}\bigg(\frac{\hat h}{m^{\hat\lambda}}\bigg)^{r-\frac{3}{2}}\big|\big|
\hat\nabla \times \hat{\bf E}\big|\big|_{r+1,\omega,\hat\Omega},
 \en
where $C_1=C_1(\omega \hat h)>0$ is independent of $p$ and $\hat{\bf E}$, but
increases as a function of the product $\omega \hat h$, and $C_1$ depends
on the shape of the $\hat\Omega_k\in \hat{\cal T}_{\hat h}$, the index $r$, the material parameters $\vartheta,
\varepsilon_r,\mu_r$ and the flux parameters.
\end{lemma}

The following auxiliary result, which states the transformation stability with respect to two mesh-dependent norms,
will play a key role in the derivation of the desired error estimates.

\begin{lemma} \label{importl}
For ${\bf E}\in {\bf V}({\cal T}_h)$, we have
\beq \label{imlemn}
\begin{split}
  \big|\big|\big|{\bf E}\big|\big|\big|_{\mathcal{F}_h} & \leq
\rho^{1/2} ~ m_{\text{mid}}^{-\frac{1}{2}}~m_{\text{min}}^{-\frac{1}{2}}~ (1 + ||\Lambda^{\frac{1}{2}}||^{-1} ) ~\big|\big|\big|\hat{\bf E}\big|\big|\big|_{\hat{\mathcal{F}}_{\hat h}}, \\
 \big|\big|\big|{\bf E}\big|\big|\big|_{\mathcal{F}_h^+} & \leq
 \rho^{1/2} ~ m_{\text{mid}}^{-\frac{1}{2}}~m_{\text{min}}^{-\frac{1}{2}}~ (1 + ||\Lambda^{\frac{1}{2}}||^{-1} ) ~ \big|\big|\big|\hat{\bf E}\big|\big|\big|_{\hat{\mathcal{F}}_{\hat h}^+}.
\end{split}
\eq
\end{lemma}

{\it Proof}.  We divide the proof into three steps.

{\it Step 1:} To estimate $\big|\big|\llbracket
\mu^{-1}\nabla_h\times {\bf E}\rrbracket_T\big|\big|_{0,\mathcal{F}_h^I}$ and $~\big|\big| {\bf n}\times(\mu^{-1}\nabla_h\times {\bf
E})\big|\big|_{0,\mathcal{F}_h^B}$.

By the first equation of (\ref{transrela1}) and the first equation of (\ref{transrela16}), we can verify that
\be \label{ne2}
\mu^{-1}\nabla_h\times{\bf E} =
\mu_r^{-1}~P^T\Lambda^{-\frac{1}{2}}~\hat\nabla_{\hat h}\times \hat{\bf E}.
\en
Thus, on the interface $\Gamma_{kj}\in \mathcal{F}_h^I$ we have
\be \label{ne3}
\llbracket \mu^{-1}\nabla_h\times {\bf E}\rrbracket_T =
\mu_r^{-1}~ \big( {\bf n}_k\times (P^T\Lambda^{-\frac{1}{2}}~\hat\nabla_{\hat h}\times \hat{\bf E}_k) + {\bf n}_j\times (P^T\Lambda^{-\frac{1}{2}}~\hat\nabla_{\hat h}\times \hat{\bf E}_j) \big).
\en
It follows by the transformation (\ref{trans}) that
\be \label{ne4}
 {\bf n}_k = |S^{-T}{\bf n}_k|~S^T\hat{\bf n}_k=|S^{-T}{\bf n}_k|~P^TM^T\hat{\bf n}_k.
\en
Notice that $(P^T{\bf a})\times(P^T{\bf b})= P^T({\bf a}\times{\bf b})~(\forall{\bf a},{\bf b} \in R^3)$. Then, by (\ref{ne4}) we obtain
\beq \label{ne5}
\begin{split}
{\bf n}_k &\times  (P^T\Lambda^{-\frac{1}{2}}~\nabla_{\hat h}\times \hat{\bf E}_k)
 = |S^{-T}{\bf n}_k|~(P^TM^T\hat{\bf n}_k) \times (P^T\Lambda^{-\frac{1}{2}}~\hat\nabla_{\hat h}\times \hat{\bf E}_k) &
 \\
 &= |S^{-T}{\bf n}_k|~P^T\bigg((M^T\hat{\bf n}_k) \times (\Lambda^{-\frac{1}{2}}~\hat\nabla_{\hat h}\times \hat{\bf E}_k)\bigg)
 = |S^{-T}{\bf n}_k|~ P^T\Lambda^{1/2} \hat{\bf n}_k \times (\hat\nabla_{\hat h}\times\hat{\bf E}_k).&
\end{split}
\eq
Substituting (\ref{ne5}) into (\ref{ne3}), yields
\be \label{ne7}
\llbracket \mu^{-1}\nabla_h\times {\bf E}\rrbracket_T =
|S^{-T}{\bf n}_k|P^T\Lambda^{1/2}
\llbracket \mu_r^{-1}(\hat\nabla_{\hat h}\times\hat{\bf E})\rrbracket_T.
\en
It is easy to see that
\be \label{inequ100}
|S^{-T}{\bf n}_k|~||\Lambda^{1/2}|| \leq \rho^{1/2},
\en
and
\be \label{ne7_7}
\bigg(\frac{|\Gamma_{kj}|}{|\hat{\Gamma}_{kj}|}\bigg)^{\frac{1}{2}}\leq ~ m_{\text{mid}}^{-\frac{1}{2}}~m_{\text{min}}^{-\frac{1}{2}}.
\en
These, together with (\ref{ne7}), lead to
\beq \label{s1_1}
\begin{split}
\big|\big|\llbracket \mu^{-1}\nabla_h\times {\bf E}\rrbracket_T \big|\big|_{0,\mathcal{F}_h^I}
 & \leq ~ \rho^{1/2} ~ m_{\text{mid}}^{-\frac{1}{2}}~m_{\text{min}}^{-\frac{1}{2}} ~
 \big|\big|\llbracket \mu_r^{-1}\hat\nabla_{\hat h}\times
\hat{\bf E}\rrbracket_T\big|\big|_{0,\hat{\mathcal{F}}_{\hat h}^I}&
\end{split}
\eq
 and
 \begin{eqnarray} \label{s1_2}
  & \quad\quad\big|\big| {\bf n}\times(\mu^{-1}\nabla_h\times {\bf
E})\big|\big|_{0,\mathcal{F}_h^B} \leq ~\rho^{1/2}~ m_{\text{mid}}^{-\frac{1}{2}}~m_{\text{min}}^{-\frac{1}{2}}~
\big|\big|\hat{\bf n}\times(\mu_r^{-1}\hat\nabla_{\hat h}\times \hat{\bf
E})\big|\big|_{0,\hat{\mathcal{F}}_{\hat h}^B}.
 \end{eqnarray}

{\it Step 2:} Build estimates of $\big|\big|\llbracket {\bf
E}\rrbracket_T \big|\big|_{0,\mathcal{F}_h^I}$ and $\big|\big| {\bf
n}\times {\bf E}\big|\big|_{0,\mathcal{F}_h^B}. $

By (\ref{ne4}), the scaling transformation (\ref{standa2}) and the coordinate transformation (\ref{trans}), we
obtain
\be \label{ne9}
{\bf n}_k\times{\bf E}_k = |S^{-T}{\bf n}_k|~(P^TM^T\hat{\bf n}_k)
\times (P^T\Lambda^{-\frac{1}{2}}\hat{\bf E}_k). \en
\textcolor{red}{By direct manipulation,} we can show that
\be \label{ne10}
 {\bf n}_k\times{\bf
E}_k =  |S^{-T}{\bf n}_k|~ P^T\Lambda^{1/2} ~\hat{\bf n}_k\times\hat{\bf E}_k.
\en
Combining (\ref{ne10}) with (\ref{inequ100}) and (\ref{ne7_7}) yields
\begin{eqnarray} \label{ne11}
& \big|\big|\llbracket {\bf E}\rrbracket_T \big|\big|_{0,\mathcal{F}_h^I} \leq  ~\rho^{1/2}~ m_{\text{mid}}^{-\frac{1}{2}}~m_{\text{min}}^{-\frac{1}{2}}~
\big|\big|\llbracket \hat{\bf E}\rrbracket_T \big|\big|_{0,\hat{\mathcal{F}}_{\hat h}^I},
\cr &  \big|\big|{\bf n}\times {\bf E}\big|\big|_{0,\mathcal{F}_h^B} \leq ~\rho^{1/2}~ m_{\text{mid}}^{-\frac{1}{2}}~m_{\text{min}}^{-\frac{1}{2}}~
\big|\big| \hat{\bf n}\times \hat{\bf E}\big|\big|_{0,\hat{\mathcal{F}}_{\hat h}^B}.
 \end{eqnarray}

{\it Step 3:} To estimate $ \big|\big| \{\{ {\bf E} \}\}
\big|\big|_{0,\mathcal{F}_h^I}$, ~$ \big|\big| \{\{
\mu^{-1}\nabla_h\times {\bf E} \}\} \big|\big|_{0,\mathcal{F}_h^I}$
and $\big|\big|  {\bf n}\times{\bf E} \big|\big|_{0,\mathcal{F}_h^B}$.

By the scaling transformation (\ref{standa2}) and the coordinate
transformation (\ref{trans}), we get
\be \label{ne12}
\big|\big| \{\{ {\bf E}
\}\} \big|\big|_{0,\mathcal{F}_h^I} =  \big|\big| \{\{ G\tilde{\bf
E} \}\} \big|\big|_{0,\mathcal{F}_h^I}  \leq ||G||~
\big|\big| \{\{ \tilde{\bf E} \}\} \big|\big|_{0,\mathcal{F}_h^I}
 ~\leq~||\Lambda^{-\frac{1}{2}}||~ m_{\text{mid}}^{-\frac{1}{2}}~m_{\text{min}}^{-\frac{1}{2}}~ \big|\big| \{\{ \hat{\bf E} \}\} \big|\big|_{0,\hat{\mathcal{F}}_{\hat h}^I}.
\en
Moreover, we have
\beq \label{ne13}
\begin{split}
\big|\big| \{\{ \mu^{-1}\nabla_h\times {\bf E} \}\}
\big|\big|_{0,\mathcal{F}_h^I} & \xlongequal{(\ref{ne2})}
 \big|\big| \{\{ \mu_r^{-1}~P^T\Lambda^{-\frac{1}{2}}~\hat\nabla_{\hat h}\times \hat{\bf E} \}\} \big|\big|_{0,\mathcal{F}_h^I}&
 \\
&\leq ~||\Lambda^{-\frac{1}{2}}||~ m_{\text{mid}}^{-\frac{1}{2}}~m_{\text{min}}^{-\frac{1}{2}}~ \big|\big| \{\{ \mu_r^{-1}\hat\nabla_{\hat h}\times \hat{\bf E} \}\} \big|\big|_{0,\hat{\mathcal{F}}_{\hat h}^I},&
\end{split}
\eq
and
\begin{eqnarray}\label{ne14}
& \big|\big|   {\bf n}\times{\bf E}  \big|\big|_{0,\mathcal{F}_h^B} \overset{(\ref{ne11})} {\leq}
~\rho^{1/2}~ m_{\text{mid}}^{-\frac{1}{2}}~m_{\text{min}}^{-\frac{1}{2}}~
\big|\big| \hat{\bf n}\times\hat{\bf E} \big|\big|_{0,\hat{\mathcal{F}}_{\hat h}^B}.
\end{eqnarray}
Combining (\ref{s1_1})-(\ref{s1_2}) with (\ref{ne11})-(\ref{ne14}) gives the desired results (\ref{imlemn}). 

 $\Box$
\begin{remark} We point out that, since we have used different choice of $G$ from that in \cite{yuan_aniso}, the transformation stability estimates (\ref{imlemn})
with respect to the condition number $\rho$ of $\Lambda$ are better than (5.7) and (5.8) in \cite{yuan_aniso}.
\end{remark}

\begin{remark}
As in the proof of the above Lemma, we can obtain the following transformation stability with respect to two mesh-dependent norms, for $\forall \hat {\bf E} \in \
\hat{\bf V}(\hat{\cal T}_{\hat h})$,
\beq \label{timlemn}
\begin{split}
  \big|\big|\big| \hat{\bf E} \big|\big|\big|_{\hat{\mathcal{F}}_{\hat h}} & \leq
\rho^{1/2} ~ m_{\text{max}}^{\frac{1}{2}}~m_{\text{mid}}^{\frac{1}{2}}~ (1 + ||\Lambda^{-\frac{1}{2}}||^{-1} ) ~ \big|\big|\big| {\bf E} \big|\big|\big|_{\mathcal{F}_h}, \\
 \big|\big|\big| \hat{\bf E} \big|\big|\big|_{\hat{\mathcal{F}}_{\hat h}^+} & \leq
 \rho^{1/2} ~ m_{\text{max}}^{\frac{1}{2}}~m_{\text{mid}}^{\frac{1}{2}}~ (1 + ||\Lambda^{-\frac{1}{2}}||^{-1} ) ~ \big|\big|\big| {\bf E} \big|\big|\big|_{\mathcal{F}_h^+}.
\end{split}
\eq
\end{remark}

\textcolor{red}{
As in \cite{pwdg}, we define the following slightly modified weaker norm than $L^2$-norm in $\hat\Omega$:  for every   $\hat{\bf v}\in L^2(\hat\Omega)^3$,
\be \label{tdivdefi}
\big|\big|\hat{\bf v}\big|\big|_{{\bf H}(\text{div};\hat\Omega)'}
=\mathop{\text{sup}~~}\limits_{\hat{\bf
u}\in {\bf H}(\text{div};\hat\Omega)} \frac{\int_{\hat\Omega}\hat{\bf v}\cdot\hat{\bf u}~d\hat{\bf
x}}{\big|\big|\hat{\bf u}\big|\big|_{{\bf H}(\text{div};\hat\Omega)}},
\en
where
\be \label{ldivnorm}
\big|\big|\hat{\bf u}\big|\big|_{{\bf H}(\text{div};\hat\Omega)}^2=\big|\big|\hat{\bf
v}\big|\big|_{0,\hat\Omega}^2 + \big|\big|\hat\nabla\cdot\hat{\bf
v}\big|\big|_{0,\hat\Omega}^2,
\en
Notice that we have $\big|\big|\hat{\bf v}\big|\big|_{{\bf H}(\text{div};\hat\Omega)'}= \big|\big| \hat{\bf v}\big|\big|_{0,\Omega}$ for every
$$ \hat{\bf v}\in {\bf H}(\text{div}^0;\hat\Omega) = \bigg\{\hat{\bf v}\in L^2(\hat\Omega)^3: \nabla\cdot \hat{\bf v}=0 ~\text{in} ~\hat\Omega \bigg\}. $$
}

Define the following norm: for every ${\bf w}$ satisfying $G^{-1}{\bf w}\in L^2(\Omega)^3$,
\be \label{divdefi}
\big|\big|{\bf
w}\big|\big|_{\tilde{\bf H}(\text{div};\Omega)'}: =
\mathop{\text{sup}~~}\limits_{{\bf v}\in
\tilde{\bf H}(\text{div};\Omega)} \frac{\int_{\Omega}G^{-1}{\bf
w}\cdot G^{-1}{\bf v}~d{\bf x}}{\big|\big|{\bf
v}\big|\big|_{\tilde{\bf H}(\text{div};\Omega)}}
\en
where
\be \label{fdivnorm}
\big|\big|{\bf v}\big|\big|_{\tilde{\bf H}(\text{div};\Omega)}^2 = \big|\big| G^{-1}  {\bf
v}\big|\big|_{0,\Omega}^2 + \big|\big|\nabla\cdot \mathcal{B}{\bf v}\big|\big|_{0,\Omega}^2
\en
with $\mathcal{B} = \frac{A} {\text{det}^{\frac{1}{2}}(A)}.$

By  direct calculation, we have
\be  \label{divrela}
\big|\big|\nabla\cdot \mathcal{B}{\bf v}\big|\big|_{0,\Omega}^2 = \big|\big|\nabla\cdot \mathcal{B}(G\hat{\bf v})\big|\big|_{0,\Omega}^2
 = \big|\big|\nabla\cdot (P^TM^{-1}\hat{\bf v})\big|\big|_{0,\Omega}^2 = \text{det}(M^{-1}) \big|\big|\hat\nabla\cdot\hat{\bf
v}\big|\big|_{0,\hat\Omega}^2,
\en
which yields
\be \label{divrelaequal}
\big|\big|{\bf v}\big|\big|_{\tilde{\bf H}(\text{div};\Omega)}^2 =
\text{det}(M^{-1}) \big|\big|\hat{\bf v}\big|\big|_{{\bf H}(\text{div};\hat\Omega)}^2.
\en

By the scaling transformation (\ref{standa2}) and the coordinate
transformation (\ref{trans}), we obtain the
following Lemma.
\begin{lemma} \label{divnorminsta}
For ${\bf E} \in {\bf V}({\cal T}_h)$, we have
 \be\label{indemeshnorm}
  \big|\big|{\bf
E}\big|\big|_{\tilde{\bf H}(\text{div};\Omega)'} =
~\text{det}(M^{-1/2})~
\big|\big|\hat{\bf E}\big|\big|_{{\bf H}(\text{div};\hat\Omega)'}.
 \en
\end{lemma}

{\it Proof}. By the scaling argument, we have
\begin{eqnarray} \nonumber
& \big|\big|{\bf E}\big|\big|_{\tilde{\bf H}(\text{div};\Omega)'}
\xlongequal{(\ref{divdefi})} \mathop{\text{sup}~~}\limits_{{\bf
v}\in \tilde{\bf H}(\text{div};\Omega)}
\frac{\int_{\Omega}G^{-1}{\bf E}\cdot G^{-1}{\bf v}~d{\bf x}}
{\big|\big|{\bf v}\big|\big|_{\tilde{\bf H}(\text{div};\Omega)}}
\xlongequal{(\ref{divrelaequal})} \mathop{\text{sup}~~}\limits_{\hat{\bf
v}\in {\bf H}(\text{div};\hat\Omega)} \text{det}(M^{-1/2})~
\frac{\int_{\hat\Omega}\hat{\bf E}\cdot\hat{\bf v}~d\hat{\bf
x}}{\big|\big|\hat{\bf v}\big|\big|_{{\bf H}(\text{div};\hat\Omega)}}.
\end{eqnarray}

 $\Box$

With the help of the above preparation, we can prove the final results easily.

\begin{theorem} \label{the1} Let ${\bf E}$ and ${\bf E}_h$ denote the analytical solution of (\ref{eq1})-(\ref{eq2}) and the solution of (\ref{maxwelldis}),
respectively. Suppose that $p$ and $r$ satisfy the conditions in Lemma \ref{approesti}. Assume \textcolor{red}{that $\omega \hat h \leq C$} and ${\bf E}\in
{\bf H}^{r+1}(\text{curl};\Omega)$. Then, for sufficiently large $p$, we have
 \be \label{fhnorm}
\big|\big|\big| {\bf E} - {\bf E}_h \big|\big|\big|_{\mathcal{F}_h}
\leq \textcolor{red}{C}~C_2 ~\rho^{\frac{3}{4}}  ~\omega^{-5/2}(\frac{
h}{m^{\hat\lambda}})^{r-\frac{3}{2}}\big|\big|\nabla \times {\bf
E}\big|\big|_{r+1,\omega,\Omega} \en
and
\be
\label{innorm2} \big|\big|{\bf E}- {\bf
E}_h\big|\big|_{\tilde{\bf H}(\text{div};\Omega)'} \leq
~\textcolor{red}{C}~C_3 ~\rho^{\frac{5}{4}}~\omega^{-3}
\frac{h^{r-2}}{m^{\hat\lambda(r-\frac{3}{2})}} \big|\big| \nabla \times {\bf E}\big|\big|_{r+1,\omega,\Omega},
\en
where
the constants $C_2$ and $C_3$ are defined as
$$ C_2 = ||M^{-1}||^{2}~ (1 + ||A^{\frac{1}{2}}||^{-1} ) ~||A^{-\frac{1}{2}}||\quad\mbox{and}\quad
C_3= ||M^{-1}||^{3}~ (1 + ||A^{\frac{1}{2}}||^{-1} ) ~(1 + ||A^{-\frac{1}{2}}|| ).$$
\end{theorem}

{\it Proof}. With the transformations (\ref{standa2}) and (\ref{trans}),
we define ${\bm \xi}_{ h}({\bf x}) = G~\hat{\bm\xi}_{\hat h}(S{\bf x})$, where $\hat{\bm\xi}_{\hat h}$ satisfying (\ref{approesti2})
denotes the plane wave approximation of the scaled electric field $\hat{\bf E}$. \textcolor{red}{Under the assumption $\omega \hat h \leq C$, we have
$C_1(\omega \hat h)\leq C$}. Thus, by (\ref{abstracterror}) and (\ref{imlemn}), together with (\ref{approesti2}), we deduce that
\begin{eqnarray} \label{inteerr3}
& \big|\big|\big| {\bf E} - {\bf E}_h \big|\big|\big|_{\mathcal{F}_h}
 \leq 3 ~\big|\big|\big| {\bf E} - {\bm \xi}_{ h} \big|\big|\big|_{{\mathcal{F}}_{ h}^+}
  \leq 3~\rho^{1/2} ~ m_{\text{mid}}^{-\frac{1}{2}}~m_{\text{min}}^{-\frac{1}{2}}~ (1 + ||\Lambda^{\frac{1}{2}}||^{-1} )~\big|\big|\big| \hat{\bf E} - \hat{\bm\xi}_{\hat h} \big|\big|\big|_{\hat{\mathcal{F}}_{\hat h}^+}
  \cr & \leq C~\rho^{1/2}  ~ m_{\text{mid}}^{-\frac{1}{2}}~m_{\text{min}}^{-\frac{1}{2}}~ (1 + ||\Lambda^{\frac{1}{2}}||^{-1} ) ~ \omega^{-5/2}(\frac{\hat h}{m^{\hat\lambda}})^{r-\frac{3}{2}}\big|\big| \hat\nabla \times \hat{\bf E}\big|\big|_{r+1,\omega,\hat\Omega}.
\end{eqnarray}
By the first equation of (\ref{transrela1}) and the first equation of (\ref{transrela16}), we can deduct that
 \be \label{he1}
  \big(\hat\nabla \times \hat{\bf E}\big)(\hat{\bf x}) =
  \Lambda^{-1/2}P~\big(\nabla \times {\bf E}\big)({\bf x}).
   \en
Substituting (\ref{relahhath}) and (\ref{he1}) into (\ref{inteerr3}), and using the scaling argument, yields
\beq \label{elecerr}
\begin{split}
\big|\big|\big| {\bf E} & - {\bf E}_h \big|\big|\big|_{\mathcal{F}_h}
 \leq  C ~\rho^{1/2}  ~ m_{\text{mid}}^{-\frac{1}{2}}~m_{\text{min}}^{-\frac{1}{2}}~ (1 + ||\Lambda^{\frac{1}{2}}||^{-1} ) ~ \omega^{-5/2}(\frac{\hat h}{m^{\hat\lambda}})^{r-\frac{3}{2}}~
 \big|\big| \Lambda^{-1/2}P~\nabla \times {\bf E}\big|\big|_{r+1,\omega,\hat\Omega} &
\\ & \leq
C ~\rho^{1/2} ~ m_{\text{mid}}^{-\frac{1}{2}}~m_{\text{min}}^{-\frac{1}{2}}~(1 + ||\Lambda^{\frac{1}{2}}||^{-1} ) ~
 \omega^{-5/2}(\frac{h}{m^{\hat\lambda}})^{r-\frac{3}{2}}~ ||M^{-1}||^{\frac{5}{2}} ~||\Lambda^{-\frac{1}{2}}|| ~  (det(M))^{1/2}  ~
\\ & \cdot \big|\big| \nabla \times {\bf E}\big|\big|_{r+1,\omega,\Omega}  ~
\leq ~C~C_2 ~\rho^{\frac{3}{4}}  ~\omega^{-5/2}(\frac{
h}{m^{\hat\lambda}})^{r-\frac{3}{2}}\big|\big|\nabla \times {\bf
E}\big|\big|_{r+1,\omega,\Omega}.&
\end{split}
\eq

Furthermore, by (\ref{indemeshnorm}), Proposition 4.8 of \cite{pwdg}, (\ref{relahhath}), (\ref{timlemn}) and (\ref{elecerr}), we obtain
\beqx \nonumber
\big|\big|{\bf E}&- &{\bf E}_h\big|\big|_{\tilde{\bf H}(\text{div};\Omega)'}
 \overset{(\ref{indemeshnorm})} {=} ~\text{det}(M^{-1/2})~ \big|\big|\hat{\bf E}- \hat{\bf E}_{\hat h} \big|\big|_{{\bf H}(\text{div};\hat\Omega)'}
\\ & \leq& C~\text{det}(M^{-1/2})~ \omega^{-\frac{1}{2}}~\hat h^{-\frac{1}{2}} \big|\big|\big| \hat{\bf E}- \hat{\bf E}_{\hat h} \big|\big|\big|_{\hat{\mathcal{F}}_{\hat h}}
\\ & \overset{(\ref{relahhath})} {\leq} & C ~\text{det}(M^{-1/2})~\omega^{-\frac{1}{2}}~h^{-\frac{1}{2}}~||M^{-1}||^{\frac{1}{2}}~ \big|\big|\big| \hat{\bf E}- \hat{\bf E}_{\hat h} \big|\big|\big|_{\hat{\mathcal{F}}_{\hat h}} 
\\ &  \overset{(\ref{timlemn})} {\leq} & C~\text{det}(M^{-1/2})~\omega^{-\frac{1}{2}}~h^{-\frac{1}{2}} ~||M^{-1}||^{\frac{1}{2}} ~ \rho^{1/2} ~ m_{\text{max}}^{\frac{1}{2}}~m_{\text{mid}}^{\frac{1}{2}}~ (1 + ||\Lambda^{-\frac{1}{2}}||^{-1} ) ~  \big|\big| {\bf E}- {\bf E}_{h} \big|\big|_{\mathcal{F}_{ h}}
\\ &   \overset{(\ref{elecerr})} {\leq} &~C~C_3 ~\rho^{\frac{5}{4}}~\omega^{-3}
\frac{h^{r-2}}{m^{\hat\lambda(r-\frac{3}{2})}} \big|\big| \nabla \times {\bf E}\big|\big|_{r+1,\omega,\Omega}.
\eqx

$\Box$

\begin{remark}
In the considered anisotropic case, the derived error estimates contain a factor depending on the condition number $\rho$ of the coefficient matrix $A$.
In contrast to the shape regularity assumption on ${\cal T}_h$ in \cite{yuan_aniso} (where $A$ is a diagonal matrix),
we assume that the transformed triangulation $\hat{{\cal T}}_{\hat{h}}$ is shape regular.
Under this assumption, we obtain the important property ${\hat h}\le C||M^{-1}||^{-1}h$ instead of ${\hat h}\le ||M||h$ (derived
in \cite{yuan_aniso}). By this property, we can build good stabilities of the transformation. Owe to such stabilities, we obtained obviously better
error estimates than those in \cite[Theorem 5.1]{yuan_aniso} in the sense that the error bounds
on the condition number $\rho$ is superior to that of \cite{yuan_aniso}.
Besides, we believe that the orders of the condition number $\rho$ in the error estimates are optimal
since the transformation stability estimates seem sharp.
\end{remark}






\section{Numerical experiments}

In this Section, we apply the proposed PWDG method to solve the acoustic wave equation and electromagnetic
wave propagation in anisotropic media, and we report some numerical
results to verify the efficiency of the proposed method.

\textcolor{red}{ For Maxwell's equations, we use ``new PWDG" to represent the
proposed method in this paper, and use ``old PWDG" to represent the method introduced in \cite{yuan_aniso}, in
which the shape regularity assumption was directly imposed on ${\cal T}_h$ (\cite{yuan_aniso} only considers the case of the diagonal matrix $A$). Here, ``old PWDG" employs the same variational formulation (\ref{maxwelldis}) and coordinate transformation (\ref{trans}) as ``new PWDG". For the Helmholtz equation, the formulation of  ``the old PWDG" is the same as (\ref{helmdisvaria}) of ``the new PWDG", and  the only difference between ``old PWDG" and ``new PWDG" is  in a different ``shape regularity" assumption on the mesh.}
As described in Section 2 and Section 3, we choose the same number \(p\) of basis functions for every elements \(\Omega_k\).
We consider the choice of numerical fluxes for the PWDG method as in  \cite{ref21}:
the constant parameters $\alpha =\beta = \delta = 1/2$.

To measure the accuracy of the numerical solutions $u_h$ and \({\bf E}_h\),
we introduce the following $L^2$ relative error:
$$\text{err.}={||u - u_h||_{L^2(\Omega)}\over{||u||_{L^2(\Omega)}}}, ~~\text{or} ~~ \text{err.}={||{\bf E}-{\bf E}_h||_{L^2(\Omega)}\over{||{\bf E}||_{L^2(\Omega)}}}  $$
for the exact solution $u\in L^2(\Omega)$, or \( {\bf E} \in (L^2(\Omega))^3\).
All of the computations have been done in MATLAB, and the system matrix was computed by exact integration on the mesh skeleton.
\textcolor{red}{The mesh ${\cal T}_{h}$ of $\Omega$ is generated by the Mesh Generation Algorithm of section 2.3.1.
We start by covering $\hat\Omega$ by a mesh of tetrahedra denoted by $\hat{{\cal T}}_{\hat{h}}$, which is generated by the software Gmsh \cite{Geuzaine}.
``DOFs" represents the number of degree of freedoms equal to the elements multiplied by the number of basis functions per element.}

\subsection{ The 3D anisotropic Helmholtz equation}


The exact solution of the problem is
$$ u({\bf x}) =
\frac{\text{exp}\bigg(\text{i}\omega\big|S~{\bf
x}-{\bf x}_0\big|\bigg)}{4\pi\big|S~{\bf x}-{\bf x}_0\big|},  $$
where \(\Omega=[0,1]^3\).
We set the anisotropic matrix $A = \left(
  \begin{array}{ccc}
  a^2+\rho^{-1}b^2 & ab(\rho^{-1}-1) & 0 \\
   ab(\rho^{-1}-1)  & b^2+\rho^{-1}a^2 & 0\\
   0 & 0 & \rho^{-1}\\
  \end{array}
\right)$, where $a=\frac{1}{\sqrt{2}},~b=(1-a^2)^{\frac{1}{2}}=\frac{1}{\sqrt{2}}$, and $\rho \geq 1$.
By orthogonal diagonalizing,  we obtain the diagonal matrix $\Lambda=\text{diag}(1,\rho^{-1},\rho^{-1})$,
and the orthogonal matrix
$P = \left(
  \begin{array}{ccc}
  a & -b & 0 \\
   b & a & 0\\
   0 & 0 & 1\\
  \end{array}
\right)$.
Thus, we have $||A||=||\Lambda||=1$, and the condition numbers of the anisotropic matrix $A=P^T\Lambda P$ and the coordinate transformation matrix $S=\Lambda^{-\frac{1}{2}}P$
 are $\rho$ and $\rho^{\frac{1}{2}}$, respectively.
 To keep the exact solution smooth in \(\Omega\),
we choose the singularity ${\bf x}_0$ as \( {\bf x}_0=(-0.6,-0.6,-0.6)\)
satisfying ${\bf x}_0\neq S{\bf x}$ where $ \forall {\bf x} \in\Omega$.

\begin{figure}[H]
\begin{center}
\begin{tabular}{c}
 \epsfxsize=0.6\textwidth\epsffile{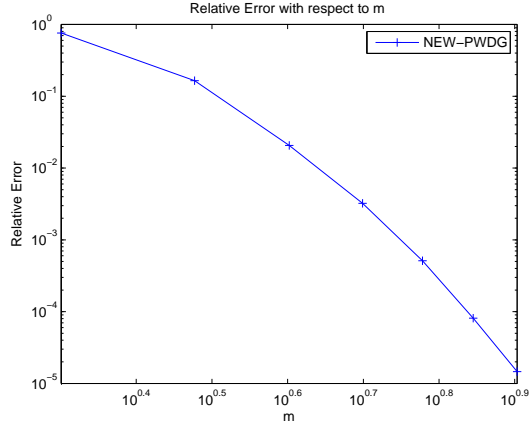}\\
\end{tabular}
\end{center}
 \caption{ \(Err.\) vs \( m \). } \label{3dhelmholtzhomog1}
\end{figure}

First, Table \ref{3dhelmholtzt1} and Figure \ref{3dhelmholtzhomog1} show the errors of the approximations
generated by the proposed PWDG method with respect to $m$. Here, we set $\omega=4\pi,~\rho=4$, choose the number $p$ of basis functions
from \(p=9\) to \(p=64\) ($2 \leq m \leq 7$), and fix the mesh triangulation.

\vskip 0.1in
\begin{center}
       \tabcaption{}
\label{3dhelmholtzt1} \vskip -0.3in
        Errors of the approximations for the case of \(p-\)convergence. \vskip 0.1in
\begin{tabular}{|c|c|c|c|c|c|c|c|} \hline
   \(p\)  & 9 & 16 & 25 & 36& 49 & 64& 81\\ \hline
 \text{err.}    & 7.61e-1 & 1.65e-1 & 2.07e-2 & 3.23e-3  & 5.12e-4  & 8.12e-5 & 1.46e-5 \\ \hline
   \end{tabular}
     \end{center}

 Figure \ref{3dhelmholtzhomog1} shows the plot of $L^2$ relative errors with
respect to \( m \). It highlights two different regimes for
increasing $m$: (i) a preasymptotic region with slow convergence,
(ii) a region of faster convergence. \textcolor{red}{The conclusion is consistent with
 Figures 4.2-4.5 of \cite{ref21}.}

\begin{figure}[H]
\begin{center}
\begin{tabular}{cc}
 \epsfxsize=0.5\textwidth\epsffile{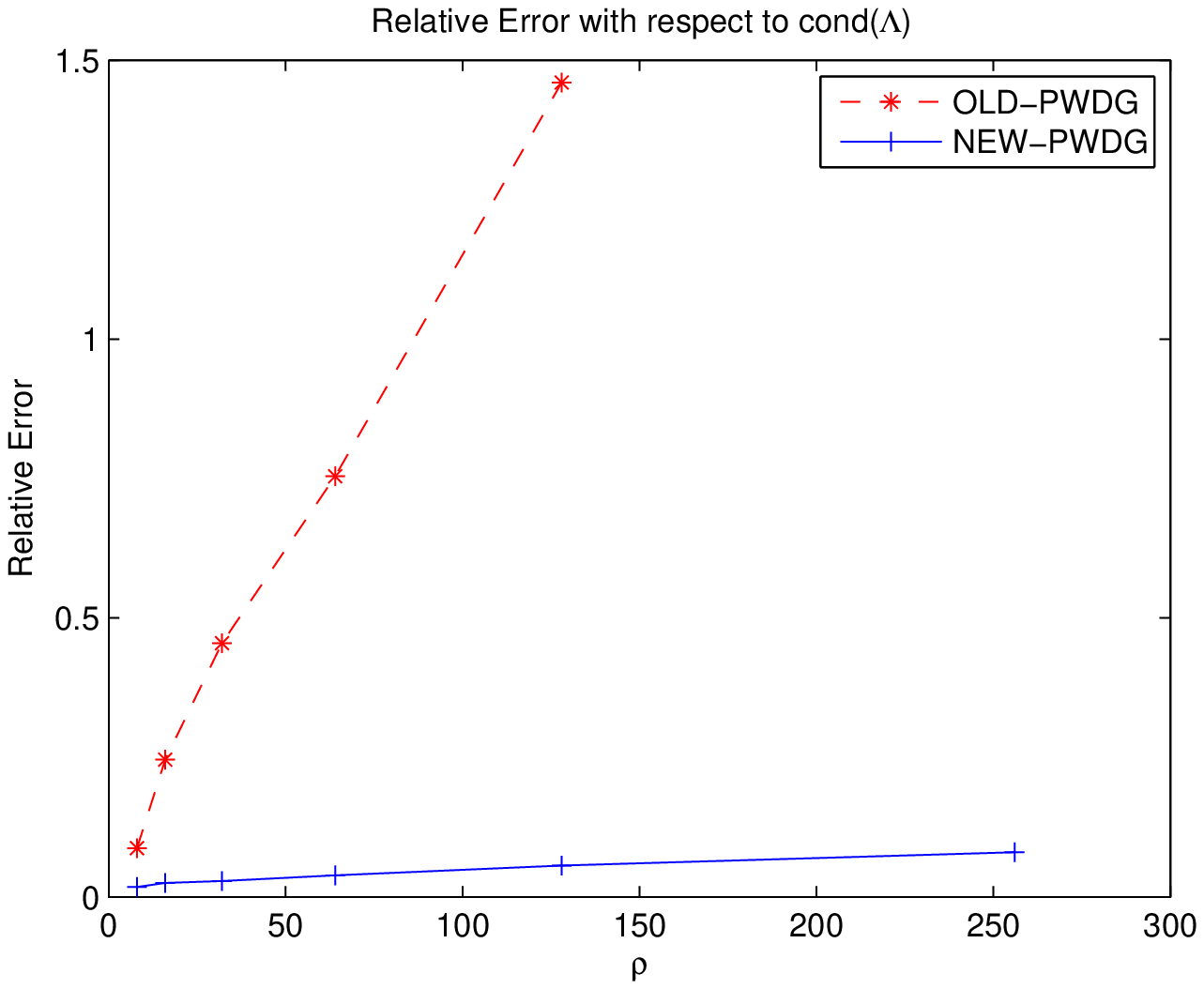}&
 \epsfxsize=0.5\textwidth\epsffile{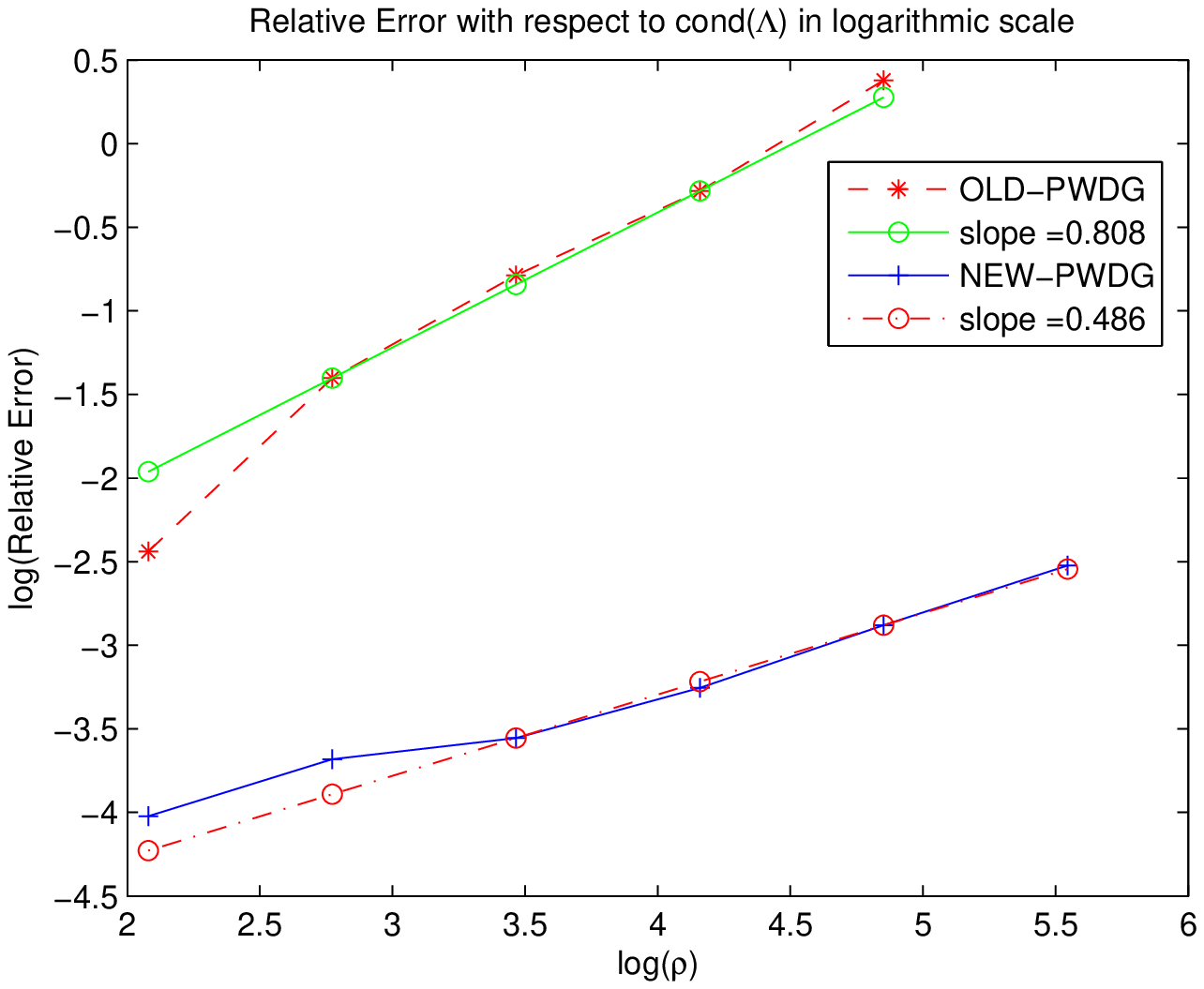}\\
\end{tabular}
\end{center}
 \caption{Left: \( Err.\) vs \(\rho\). \quad
Right: \( Err.\) vs \(\rho\) in logarithmic scale.} \label{3dhelm_rho}
\end{figure}

Then, we would like to compare the $L^2$ relative errors of the approximations generated by
the proposed PWDG method and the old PWDG method in \cite{yuan_aniso} for the case of almost the same degree of freedoms (DOFs).
We fix $\omega=4\pi, ~p=25, ~ \hat h=\frac{1}{4}$,
 but increase the condition number of the anisotropic matrix $A$ and the corresponding number of elements.
The numerical results are shown in Table \ref{HELMhomopwdgl2rhof} and Figure \ref{3dhelm_rho} .

\vskip 0.1in
\begin{center}
       \tabcaption{} \vskip -0.3in
\label{HELMhomopwdgl2rhof}
       Comparisons of errors of approximations with respect to $\rho$. \vskip 0.1in
\begin{tabular}{|c|c|c|c|c|c|c|c|} \hline
 &  \(\rho\) & 8 & 16 & 32 &  64 & 128 & 256 \\ \hline
 \multirow{2}*{ old PWDG} & err   & 8.72e-2 & 2.46e-1  & 4.55e-1  & 7.54e-1  & 1.46   & --  \\  \cline{2-8}
   & DOFs  & 12800  & 25000  & 54925  & 102400  &  200000  & 390625  \\ \hline
  \multirow{2}*{ new PWDG} & err  & 1.79e-2  & 2.52e-2  & 2.86e-2  & 3.86e-2  & 5.61e-2   & 8.02e-2  \\  \cline{2-8}
   & DOFs  &  14400 & 25600  & 52900  & 102400  &  193600  & 409600 \\ \hline
   \end{tabular}
     \end{center}

We can see from Table \ref{HELMhomopwdgl2rhof} and Figure \ref{3dhelm_rho} that,
the approximations generated by the new PWDG are more accurate than those generated by the old PWDG.
Besides, \textcolor{red}{the growth rates with respect to $\rho$ }of the $L^2$ relative errors of the approximations generated by the new PWDG method
 are much smaller than those generated by the old PWDG method, which verifies the validity of Remark \ref{helml2comp}.
In addition, the numerical convergence order of the approximation $u_h$ generated by the proposed method with respect to $\rho$
is superior to the theoretical convergence order.

\begin{figure}[H]
\begin{center}
\begin{tabular}{c}
\epsfxsize=0.6\textwidth\epsffile{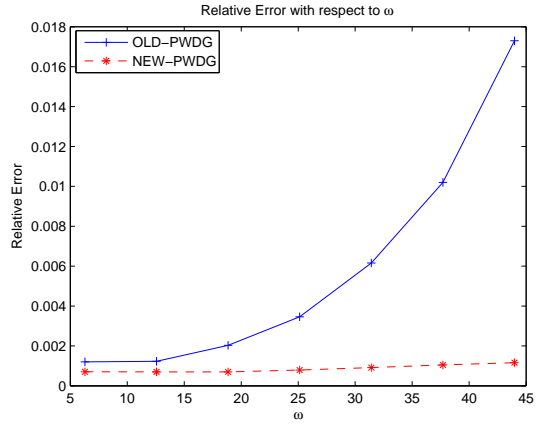}\\
\end{tabular}
\end{center}
 \caption{\( Err.\) vs \(\omega\).} \label{3dHELMhomog3_1}
\end{figure}

\vskip 0.1in
Furthermore, we fix the product of $\omega \hat h$ to be $\frac{\pi}{2}$, but increase
the wave number $\omega$ and decrease the mesh size $ \hat h$. The resulting errors of the approximations
generated by the new PWDG and the old PWDG are listed in Table \ref{HELMPWDGCOM2}
and Figure \ref{3dHELMhomog3_1}.

\vskip 0.1in
\begin{center}
       \tabcaption{}
\label{HELMPWDGCOM2} \vskip -0.3in
       Comparisons of errors of approximations with respect to $\omega$.
\vskip 0.1in
\begin{tabular}{|c|c|c|c|c|c|c|c|} \hline
 & \(\omega\)  &  $4\pi$  &  $6\pi$ & $8\pi$ & $10\pi$ & $12\pi$ & $14\pi$ \\ \hline
   \multirow{2}*{\text{old PWDG}} & err   & 1.23e-3  & 2.03e-3   & 3.46e-3 & 6.16e-3 & 1.02e-2 & 1.73e-2
     \\  \cline{2-8}
  & DOFs  & 25000 & 84375 & 200000 & 390625 & 675000 & 1071875  \\ \hline
   \multirow{2}*{\text{new PWDG}} & err  & 6.97e-4 & 6.99e-4   & 7.94e-4  & 9.13e-4  & 1.05e-3 & 1.16e-3 \\ \cline{2-8}
  & DOFs  & 25600   & 86400 & 204800 & 400000 & 691200 & 1097600 \\ \hline
   \end{tabular}
     \end{center}

Table \ref{HELMPWDGCOM2} and Figure \ref{3dHELMhomog3_1} show that, the approximations generated
by the new PWDG are more accurate than those generated by the old PWDG.
Besides, the numerical errors in relative $L^2$ norm indicate that the proposed PWDG method is slightly affected by the
pollution effect (see \cite{zhu}).

\subsection{Electric dipole in free space for a smooth case}

We compute the electric field due to an electric dipole source at
the point \( {\bf x}_0=(-0.6,-0.6,-0.6)\). The dipole point source
can be defined as the solution of a homogeneous Maxwell system
(\ref{eq1}). The exact solution of the problem is
\begin{equation}
{\bf E}_{\text{ex}}=-\text{i}\omega I \phi({\bf x},{\bf x}_0)G~{\bf
a} +\frac{I}{{\bf i}\omega\varepsilon_r} G~\nabla_{\hat h}(\nabla_{\hat h}\phi\cdot{\bf a})
\end{equation}
where
$$ \phi({\bf x},{\bf x}_0) =
\frac{\text{exp}\bigg(\text{i}\omega\sqrt{\varepsilon_r}\big|S~{\bf
x}-{\bf x}_0\big|\bigg)}{4\pi\big|S~{\bf x}-{\bf x}_0\big|}  $$ and
\(\Omega=[0,1]^3\). We set the anisotropic matrix $A = \left(
  \begin{array}{ccc}
  a^2+\rho b^2 & ab(\rho-1) & 0 \\
   ab(\rho-1)  & b^2+\rho a^2 & 0\\
   0 & 0 & \rho\\
  \end{array}
\right)$ and \(\varepsilon_r=\mu_r=1\), where $a,~b,~\rho$ are the same as in the previous Subsection.
By orthogonal diagonalizing, the diagonal matrix $\Lambda=\text{diag}(1,\rho,\rho)$, and the orthogonal matrix $P$ is the same as in the previous Section.
Besides, we have $M=\text{diag}(\rho,~\rho^{\frac{1}{2}},~\rho^{\frac{1}{2}})$ and $||M^{-1}||=\rho^{-\frac{1}{2}}$.
In addition, the condition numbers of the anisotropic matrix $A=P^T\Lambda P$ and the coordinate transformation matrix $S=MP$
 are $\rho$ and $\rho^{\frac{1}{2}}$, respectively.

\begin{figure}[H]
\begin{center}
\begin{tabular}{c}
 \epsfxsize=0.6\textwidth\epsffile{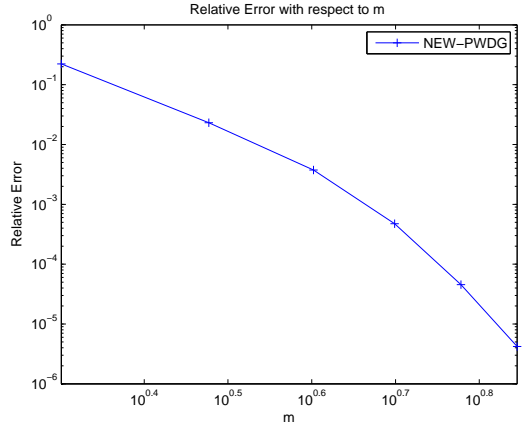}\\
\end{tabular}
\end{center}
 \caption{ \(Err.\) vs \( m \).  } \label{3dhomog1}
\end{figure}

First, Table \ref{3dmaxt1} and Figure \ref{3dhomog1} show the errors of the approximations
generated by the PWDG method with respect to $m$. Here, we set $\omega=4\pi, ~\rho=4$.
 The number $p$ of basis functions is chosen from \(p=9\) to \(p=64\) ($2 \leq m \leq 7$),
and the mesh triangulation is fixed.

\vskip 0.1in
\begin{center}
       \tabcaption{}
\label{3dmaxt1} \vskip -0.3in
        Errors of the approximations for the case of \(p-\)convergence. \vskip 0.1in
\begin{tabular}{|c|c|c|c|c|c|c|} \hline
   \(p\)  & 9 & 16 & 25 & 36& 49 & 64 \\ \hline
 \text{err.}    & 2.22e-1 & 2.31e-2 & 3.74e-3 & 4.74e-4  &  4.54e-5 & 4.20e-6 \\ \hline
   \end{tabular}
     \end{center}

Figure \ref{3dhomog1} shows the plot of $L^2$ relative errors with
respect to \( m \). It also highlights two different regimes for
increasing $m$: (i) a preasymptotic region with slow convergence,
(ii) a region of faster convergence. As stated in Remark 3.14 of \cite{ref21}, the convergence order of the approximations with
respect to $m$ turns out to be exponential since the analytical
solution of the problem can be extended analytically outside the
domain.

\begin{figure}[H]
\begin{center}
\begin{tabular}{cc}
 \epsfxsize=0.5\textwidth\epsffile{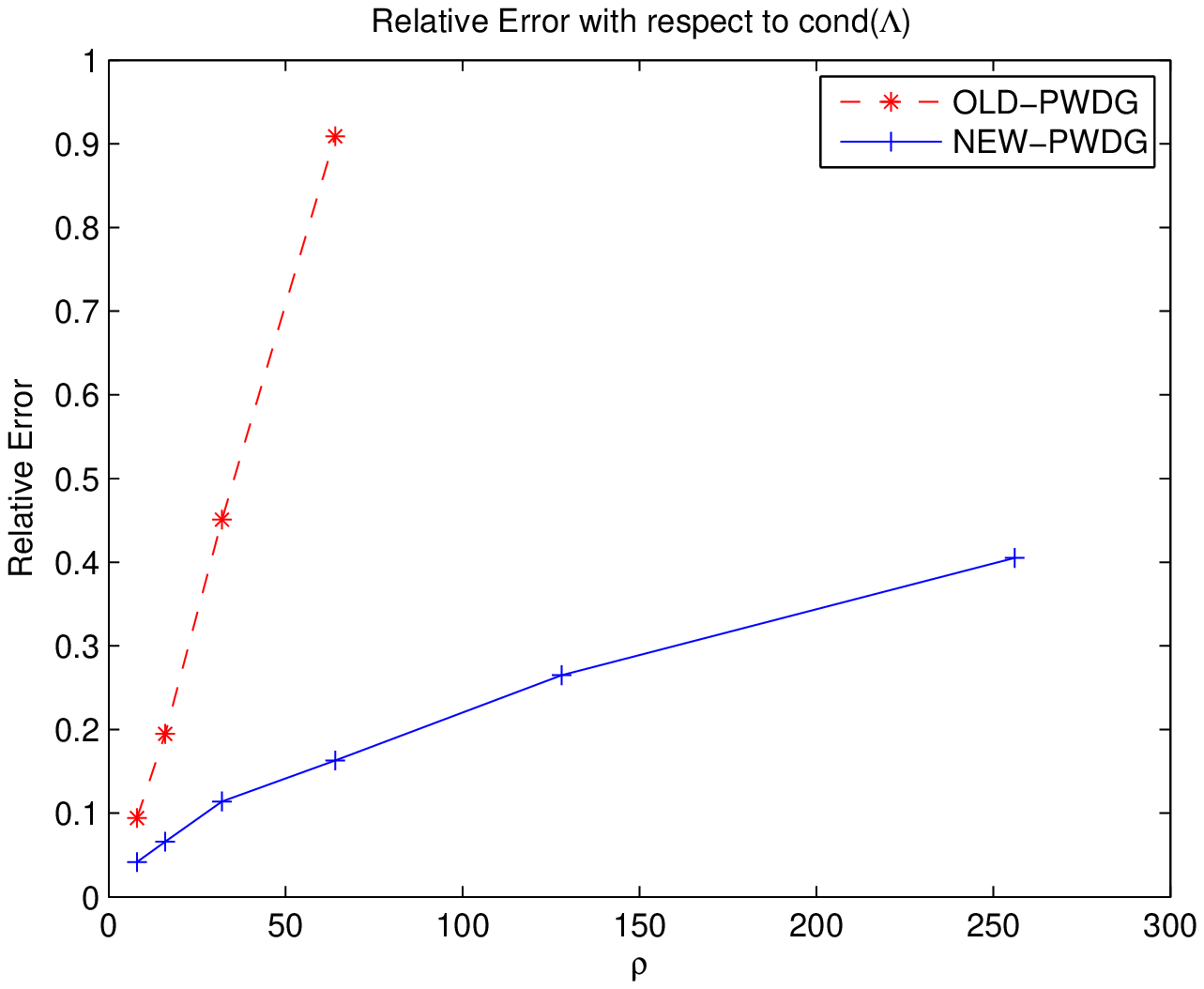}&
 \epsfxsize=0.5\textwidth\epsffile{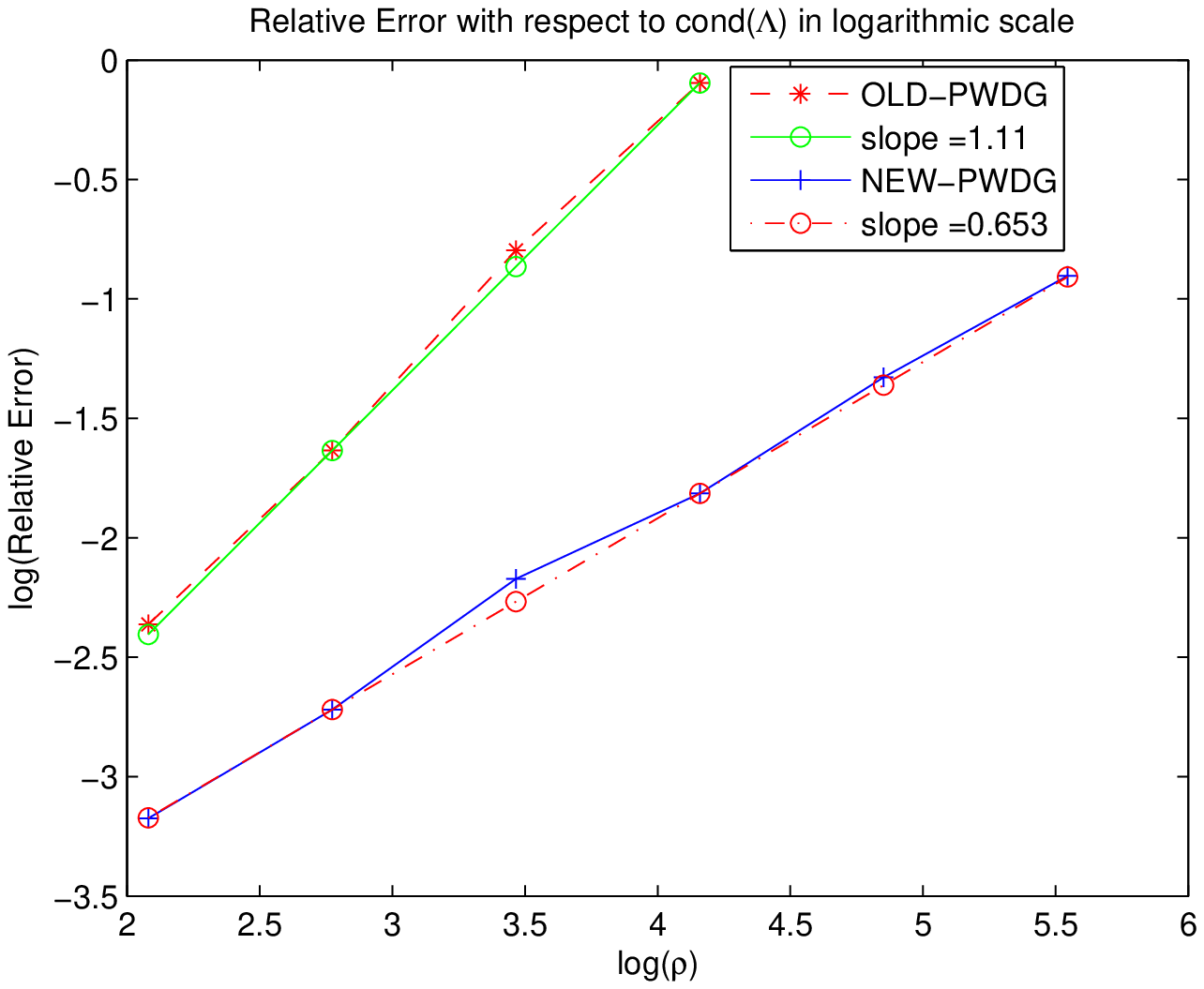}\\
\end{tabular}
\end{center}
 \caption{Left: \( Err.\) vs \(\rho\). \quad
Right: \( Err.\) vs \(\rho\) in logarithmic scale.} \label{3dmax_rho}
\end{figure}

Next, the $L^2$ relative errors of the approximations generated by the proposed PWDG method and the old PWDG method in \cite{yuan_aniso} are compared
for the case of almost the same DOFs. We fix $\omega=2\pi, ~p=25, ~\hat h=\frac{1}{2}$,
  but increase the condition number of the anisotropic matrix $A$ and the corresponding number of elements.
 The numerical results are shown in Table \ref{2homopwdgl2rho} and Figure \ref{3dmax_rho}.

\vskip 0.1in
\begin{center}
       \tabcaption{} \vskip -0.3in
\label{2homopwdgl2rho}
       Comparisons of errors of approximations with respect to $\rho$. \vskip 0.1in
\begin{tabular}{|c|c|c|c|c|c|c|} \hline
   &\(\rho\) & 8 & 16 & 32 &  64 & 128  \\ \hline
 \multirow{2}*{ old PWDG} & err   & 9.42e-2   & 1.95e-1 & 4.51e-1  & 9.09e-1  &  --     \\  \cline{2-7}
   & DOFs  & 25600   & 109850  &    400000  &  1638400  &   6632550    \\ \hline
  \multirow{2}*{ new PWDG} & err  &  4.18e-2 & 6.59e-2 &  1.14e-1 & 1.63e-1  & 2.65e-1    \\  \cline{2-7}
   & DOFs  & 28800  & 102400 & 387200  & 1638400  &  6771200    \\ \hline
   \end{tabular}
     \end{center}

From Table \ref{2homopwdgl2rho} and Figure \ref{3dmax_rho}, we can see that, the approximations generated by the new PWDG are more accurate than those generated by the old PWDG.
Moreover, \textcolor{red}{the growth rates with respect to $\rho$} of the $L^2$ relative errors of the approximations generated by the new PWDG method are much smaller than those generated by the old PWDG method.

\begin{figure}[H]
\begin{center}
\begin{tabular}{c}
 \epsfxsize=0.6\textwidth\epsffile{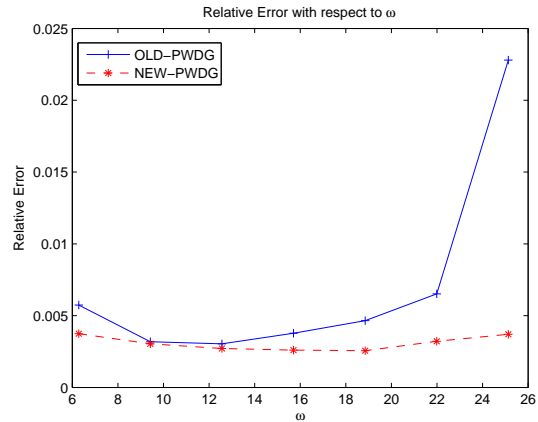}\\
\end{tabular}
\end{center}
 \caption{\( Err.\) vs \(\omega\).} \label{3dhomog3_1}
\end{figure}

Finally, we fix the product of $\omega \hat h$ to be $\frac{\pi}{2}$, but increase
the wave number $\omega$ and decrease the mesh size $\hat h$. The resulting errors of the approximations
generated by the new PWDG and the old PWDG are listed in Table \ref{PWDGCOM}
and Figure \ref{3dhomog3_1}.

\vskip 0.1in
\begin{center}
       \tabcaption{}
\label{PWDGCOM} \vskip -0.3in
       Comparisons of errors of approximations with respect to $\omega$.
\vskip 0.1in
\begin{tabular}{|c|c|c|c|c|c|c|c|} \hline
 & \(\omega\)  &  $3\pi$  &  $4\pi$ & $5\pi$ & $6\pi$ & $7\pi$ & $8\pi$ \\ \hline
      \multirow{2}*{\text{old PWDG}} & err    & 3.18e-3 &  3.04e-3 & 3.78e-3 &  4.65e-3 & 6.52e-3 & 2.28e-2  \\ \cline{2-8} 
      & DOFs  &  50000 & 109850 & 204800  & 342950  &  532400 & 781250  \\ \hline
   \text{new PWDG}  & err & 3.04e-3 &   2.71e-3  & 2.60e-3 &  2.56e-3 & 3.22e-3 & 3.70e-3  \\ \cline{2-8} 
      & DOFs  &  43200 & 102400  & 200000  &  345600 &  548800 & 819200  \\ \hline
   \end{tabular}
     \end{center}

Table \ref{PWDGCOM} and Figure \ref{3dhomog3_1} show that, the approximations generated by the new PWDG are more accurate than
those generated by the old PWDG. Besides, the numerical errors in relative $L^2$ norm indicate that the proposed PWDG method is slightly affected by the
pollution effect.

\vskip 12pt

{\bf Acknowledgments}. The authors would like to thank the anonymous reviewer, who gives many insightful comments to improve the presentation of this paper.



\end{document}